\documentclass[10pt]{amsart}
\usepackage{rotating}
\usepackage{amssymb}
\usepackage{url}
\newcommand{\q}{{\tt "}}
\newcommand{\C}{{\mathbb C}}
\newcommand{\Z}{{\mathbb Z}}
\newcommand{\cK}{\mathcal{K}}

\newcommand{\R}{{\mathbb R}}
\newcommand{\F}{{\mathbb F}}
\newcommand{\Q}{{\mathbb Q}}

\newcommand{\bu}{\bullet} 

\newcommand{\Gal}{\mbox{Gal}}

\newcommand{\ord}{\textrm{ord}}
\newcommand{\grd}{\textrm{grd}}
\newcommand{\GRD}{\textrm{grd}}
\newcommand{\rd}{\textrm{rd}}

\newcommand{\cmmt}[1]{}
\numberwithin{equation}{section}
\numberwithin{table}{section}
\numberwithin{figure}{section}
\thanks{This work was partially supported by a grant from the Simons Foundation (\#209472 to David Roberts).}

\title{A Database of Number Fields}
 \author{John W.\ Jones}
\address{School of Mathematical and Statistical Sciences, Arizona
  State University, PO Box 871804, Tempe, AZ 85287} \email{jj@asu.edu}

\author{David P.\ Roberts}
\address{Division of Science and Mathematics, University of
  Minnesota Morris, Morris, MN 56267}
\email{roberts@morris.umn.edu}
\begin{document}
\maketitle
\begin{abstract}
We describe an online database of number fields which accompanies this
paper.  The database centers 
on complete lists of number fields with prescribed invariants.  
Our description here focuses on summarizing tables and connections to
theoretical issues of current interest.
\end{abstract}

\section{Introduction}
     A natural computational problem is to completely determine the set $\cK(G,D)$ of all degree $n$
number fields $K$ with a given Galois group $G \subseteq S_n$
and a given discriminant $D$.    Many papers have 
solved instances of this problem, some relatively early contributions being
 \cite{hunter,pohst,BMO,spd}.

     This paper describes our online database of number fields at 
\begin{center}
{\tt http://hobbes.la.asu.edu/NFDB/} \; .
\end{center}
This database gives many complete determinations of $\cK(G,D)$ in
small degrees $n$, collecting 
previous results and going well beyond them. Our database complements
the  Kl\"uners-Malle online database \cite{kluners-malle}, 
which covers more groups and signatures, but is not as focused on completeness
results and the behavior of primes.   Like the Kl\"uners-Malle
database, our database is searchable and 
intralinked.  

Section~\ref{using} explains in practical terms how one can use
the database.  Section~\ref{internal} 
explains some of the internal workings of the database, including how it
keeps track of completeness.  Section~\ref{summarizing} presents tables summarizing 
the contents of the database in degrees $n \leq 11$, which is the 
setting of most of our completeness results.    
The section also briefly indicates how fields are chosen for inclusion
on the database and describes connections with 
previous work.

     The remaining sections each summarize an aspect of the 
 database, and explain how the tabulated fields  shed some light on theoretical 
 issues of current interest. 
 As a matter of
 terminology, we incorporate the signature of a field into our notion
 of discriminant, 
 considering the formal product $D = -^s |D|$ to be the
 discriminant of
 a field with $s$ complex places and absolute
 discriminant $|D|$.

   Section~\ref{quintics2357}
 focuses on the complete list of all $11279$ quintic fields 
 with Galois group $G=S_5$ and discriminant of the 
 form $-^s 2^a 3^b 5^c 7^d$.   The summarizing
 table here shows that the distribution of 
 discriminants conforms moderately well 
 to the mass heuristic of \cite{Bh}.   Section~\ref{nonsolvable567}
 summarizes lists of fields for more nonsolvable groups, but
 now with attention restricted to discriminants of the 
 form $-^s p^a q ^b$ with $p<q$ primes.  
 
   Sections~\ref{octictame} and \ref{octicwild} 
 continue to pursue cases with $D =  -^s p^a q ^b$, 
 but now for octic groups $G$ of  $2$-power order.  
 Section~\ref{octictame} treats the cases $p>2$
 and discusses connections
 to tame maximal nilpotent extensions
 as studied in \cite{boston-ellenberg, boston-perry}.  Section~\ref{octicwild} treats
 the case $p=2$ and takes a first step towards
 understanding wild ramification in some
 of the nilpotent 
 extensions studied in \cite{koch}.
 
     Sections~\ref{minnonsolvable} and \ref{gennonsolvable}
 illustrate progress in the database on  
 a  large project  initiated in \cite{jrlowgrd}.   The project
 is to completely classify Galois number fields with root discriminant 
 $|D|^{1/n}$ at most the  Serre-Odlyzko constant $\Omega := 8 \pi e^{\gamma} \approx 44.76$.
 Upper bounds on degrees coming from analysis 
 of Dedekind zeta functions \cite{martinet,od-survey} play a prominent role.  
 The database gives many solvable fields satisfying
 the root discriminant bound.  In this paper for brevity we 
 restrict attention to nonsolvable fields, where, among other interesting things, 
 modular forms \cite{bosman}, \cite{schaeffer} 
 sometimes point the way to explicit polynomials.  
   
   The database we are presenting here has its origin in posted versions
of the complete tables of our earlier work \cite{jr1}.   Other complete lists of 
fields were posted sporadically in the next ten years, while most fields 
and the new interface are  recent additions.  
Results from the predecessors of the present database
have occasionally been used as ingredients of
formal arguments, as in  e.g.\ \cite{HKK,ono-taguchi,dahmen}.  The more 
common use of our computational results has
been to guide investigations
into number fields in a more
general way.   With our recent enhancements
and this accompanying paper, we aim to 
increase the usefulness of our work to the
mathematical community.

\section{Using the database}
\label{using}
    A simple way to use the database is to request $\cK(G,D)$, for a particular
$(G,D)$.  A related but more common way is to request the union of these sets for varying 
$G$ and/or $D$.  Implicit throughout this paper and the database is that fields are always considered up
to isomorphism. As a very simple example, asking for quartic fields with any Galois group
 $G$ and discriminant $D$ satisfying  $|D| \leq 250$ returns Table~\ref{sampletable}.   

\begin{table}[htb]
\[
\begin{array}{|r r l c c l |}
\hline
\multicolumn{6}{|c|}{\mbox{\bf Results below are proven complete}} \\
\hline
\mbox{rd}(K) & \mbox{grd}(K) & D & h &  G & \mbox{Polynomial} \\
\hline
 3.29 & 6.24 & -^2 3^2 13^1 &1& D_4 &  x^4-x^3-x^2+x+1 \\
 3.34 & 3.34 & -^2 5^3 &1& C_4 & x^4-x^3+x^2-x+1 \\
 3.46 & 3.46 & -^2 2^4 3^2 &1& V_4 & x^4-x^2+1 \\
3.71 & 6.03 & -^2 3^3 7^1 &1& D_4 &  x^4-x^3+2 x+1 \\
 3.87 & 3.87 & -^2 3^2 5^2 &1& V_4 & x^4-x^3+2 x^2+x+1 \\
3.89 & 15.13 & -^2 229^1 &1&  S_4 & x^4-x+1 \\
\hline
\end{array}
\]
\caption{\label{sampletable} Results of a query for quartic fields with absolute discriminant $\leq 250$,
sorted by root discriminant} 
\end{table}
In general, the monic polynomial $f(x) \in \Z[x]$ in the last column defines the
field of its line, via $K = \Q[x]/f(x)$.  It is standardized by requiring the sum 
of the absolute squares of its complex roots to be minimal, with ties 
broken according to the conventions of Pari's {\em polredabs}.    
Note that the database, like its local analog \cite{jr-local-database}, is 
organized around non-Galois fields.  However, on a given line, 
some of the information refers to  a Galois closure $K^g$.  

The Galois group $G = \mbox{Gal}(K^g/\Q)$ is given by its common name, like in Table~\ref{sampletable},
or its $T$-name  as in \cite{butler-mckay,gp,GAP} if it does not have a very widely accepted common
name.  Information about the group---essential for intelligibility in
higher degrees---is obtainable by clicking on the group.  For example, the database reports 
$10T42$ as having structure $A_5^2.4$, hence order $60^2 4 = 14400$; moreover, it is isomorphic
to $12T278$, $20T457$, and $20T461$.  

Continuing to explain Table~\ref{sampletable}, the column $D$ prints $-^s |D|$, where $s$ is the number of
 complex
places and $|D|$ is given in factored form.   This format treats the infinite completion $\Q_\infty = \R$ on 
a parallel footing with the $p$-adic completions $\Q_p$.  If $n \leq 11$, then
clicking on any appearing prime $p$ links into 
the local database of \cite{jr-local-database}, thereby giving a detailed description 
of the $p$-adic algebra $K_p = \Q_p[x]/f(x)$.   This 
automatic $p$-adic analysis also often works in
degrees $n>11$.

The root discriminant $\mbox{rd}(K) = |D|^{1/n}$ is placed in the first column, since 
one commonly wants to sort by root discriminant.   Here and later we
often round real numbers
to the nearest hundredth without further comment.  
When it is implemented, our complete analysis at all 
ramifying primes $p$ automatically
determines the Galois root discriminant of $K$, meaning the
root discriminant of a Galois closure $K^g$.   The second
column gives this more subtle invariant $\mbox{grd}(K)$. Clicking on the entry gives
the exact form and its source.    Often it is better to sort by this 
column, as fields with the same Galois closure are then put next to 
each other.    As an example, quartic fields with $G=D_4$ come in twin pairs
with the same Galois closure.  
The twin $K^t$ of the first listed field 
$K$ on Table~\ref{sampletable} is off the table because $|D(K^t)| = 3^1 13^2 = 504$; 
however $\mbox{grd}(K^t) = \mbox{grd}(K) = 3^{1/2} 13^{1/2} \approx 6.24$.  

Class numbers are given in the column $h$, factored as $h_1 \cdots h_d$ 
where the class group is a product of cyclic groups of size $h_i$.   There is
a toggle button, so that one can alternatively receive narrow class numbers
in the same format.   To speed up the construction of the table, class
numbers were computing assuming the generalized Riemann hypothesis;
they constitute the only  part of the database that is conditional. 
Standard theoretical facts about class numbers can be 
seen repeatedly in various parts of the database.  For example, 
the unique $D_7$ septic field with Galois root discriminant 
$1837^{1/2} \approx 42.86$ has
class number $13$.  The degree $14$ Galois closure
is then forced to have class number $13 \cdot 13$ and this fact is
explicitly confirmed by the database.

When the response to a query is known to be complete, then the
table is headed by the completeness statement shown in Table~\ref{sampletable}.
   As emphasized in the introduction, keeping track
of completeness is one of the most important features of the 
database.    The completeness statement often 
reflects a very long computational proof, even if 
the table returned is very short.  

There are many other ways to search the database, mostly connected to the
behavior of primes.  For example, one can restrict the search to find
fields with restrictions on $\ord_p(D)$,  or one can search
directly for fields with Galois root discriminant in a given range.   On the
other hand, there are some standard invariants of fields that
the database does not return, such as Frobenius partitions 
and regulators.  The database does allow users to download the
list of polynomials returned, so that it can be used
as a starting point for further investigation.

\section{Internal structure}
\label{internal}

The website needs to be able to search and access a large amount of
information.  It uses a fairly standard architecture: data is stored
in a MySQL database and web pages are generated by programs written in
Perl.

A MySQL database consists of a collection of tables where each table
is analogous to a single spreadsheet with columns representing the
types of data being stored.  We use data types for integers,
floating point numbers, and strings, all of which come in various
sizes, i.e., amount of memory devoted to a single entry.
 When searching, one can
use equalities and inequalities where strings are ordered
lexicographically.

When a user requests number fields, the Perl program takes the
following steps:
\begin{enumerate}
\item Construct and execute a MySQL query to pull fields from the
  database.
\item Filter out fields which satisfy all of the user's requirements
 when needed (see below).
\item Check completeness results known to the database.
\item Generate the output web page.
\end{enumerate}

The main MySQL table has one row for each field.  There are columns
for each piece of information indicated by the input boxes in the top
portion of the search screen, plus columns for the defining polynomial (as
a string), and an internal identifier for the field.
The only unusual aspect of this portion of the database is how
discriminants are stored and searched.  The difficulty stems from the
fact that many number fields in the database have discriminants which are
too large to store in MySQL as integers.  An option would
be to store the discriminants as strings, but then it would be
difficult to search for ranges: string comparisons in MySQL are
lexicographic, so `11' comes before `4'.  Our solution is to store
absolute discriminants $|D|$ as strings, but prepend the string with four
digits which give $\lfloor\log_{10}(|D|)\rfloor$, padded on the left
with zeros as needed.  So, $4$ is stored as `00004', $11$ is `000111',
etc.  This way we can use strings to store each discriminant in
its entirety, but searches for ranges work correctly.

The MySQL table of number fields also has a column for the list of all
primes which ramify in the field, stored as a string with a separator
between primes.  This is used to accelerate searches when it is clear
from the search criteria that only a small finite list of
possibilities can occur, for example, when the user has checked the
box that ``Only listed primes can ramify''.

Information on ramification of specific primes can be input in the
bottom half of the search inputs.  To aid in searches involving these
inputs, we have a second MySQL table, the ramification table, which
stores a list of triples.  A triple $(\text{field identifier}, p, e)$ 
indicates that $p^e$ exactly divides the discriminant of the
corresponding field.  The most common inputs to the bottom half of the
search page work well with this table, namely those which list
specific primes and allowable discriminant exponents.  However, the
search boxes allow much more general inputs, i.e., where a range of
values is allowed for the prime and the discriminant exponent allows
both $0$ and positive values.  It is possible to construct MySQL
queries for inputs of this sort, but they are complicated, involve
subqueries, and are relatively slow.  Moreover, a search condition of this
type typically rules out relatively few number fields.  If a user does
make such a query, we do not use the information at this stage.
Instead, we invoke Step~(2) above to select fields from the MySQL
query which satisfy these additional requirements.

The database
supports a variety of different types of completeness results.
Complicating matters is that these results can be
 interrelated.  We use four MySQL tables for
storing ways in which the data is complete.  In describing them, $G$
denotes the Galois group of a field, $n$ is the degree, $s$ is the
number of complex places, and $|D|$ is the absolute discriminant,
as above.  The tables are
\begin{enumerate}
\item[A.] store $(n, s, B)$ to indicate that the database is complete
  for fields with the given $n$ and $s$ such that  $|D|\leq B$;
\item[B.] store $(n, s, G, B)$  to indicate that the database is complete
  for fields with the given $n$, $s$, and $G$ such that  $|D|\leq B$;
\item[C.]  store $(n, S, L)$ where $S$ is a list of primes and $L$ is a list
  of Galois groups to indicate that the database is complete for
  degree $n$ fields unramified outside $S$ for each Galois group in $L$;
\item[D.] store $(n, G, B)$ to indicate that the database is complete for
  degree $n$ fields $K$ with Galois group $G$ such that $\grd(K)\leq B$.
\end{enumerate}
In each case, database entries include the degree, so individual
Galois groups can be stored by their $T$-number (a small integer).  In
the third case, we store the list $L$ by an integer whose bits
indicate which $T$-numbers are included in the set.  For example,
there are $50$ $T$-numbers in degree $8$, so a list of Galois groups
in that degree is a subset of $S\subseteq \{1,\ldots, 50\}$ which we
represent by the integer $\sum_{t\in S} 2^{t-1}$.  These integers are
too large to store in the database as integers, so they are stored as
strings, and converted to multiprecision integers in Perl.  
The list of primes in the third table is simply stored as a string
consisting of the primes and separating characters.

To start checking for completeness, we first check that there are only
finitely many degrees involved, and that the search request contains
an upper bound on at least one of: $|D|$, $\rd(K)$, $\grd(K)$, or the
largest ramifying prime.
We then loop over the degrees in the user's search.  We allow for the
possibility that a search is known to be complete by some combination
of completeness criteria.  So throughout the check, we maintain a list
of Galois groups which need to be checked,
and the discriminant values to
check.  If one check shows that some of the Galois groups for the
search are known to be complete, they are removed from the list.  If
that list drops to being empty, then the search in that degree is
known to be complete.  Discriminant values are treated analogously.

For each degree, bounds on $|D|$ and $\rd(K)$ are clearly equivalent.
Less obviously, bounds between $\rd(K)$ and $\grd(K)$ are related.  In
particular, we always have $\rd(K)\leq \grd(K)$, but also have for
each Galois group, $\grd(K)\leq \rd(K)^{\alpha(G)}$ where $\alpha(G)$
is a rational number depending only on $G$ (see \cite{jr-tame-wild}).

We then perform the following checks.
\begin{itemize}
\item We compare the request with Tables~A, B, and D for
  discriminant bound restrictions.
\item Remove Galois groups from the list to be checked based on
  $\grd$.
\item If there are at most ten discriminants not accounted for, check
  each individually against Table~C.
\item If there is a bound on the set of ramifying primes, which could
  arise from the user checking ``Only these primes ramify'', or from a
  bound on the maximum ramifying prime, check Table~C.
\end{itemize}

\section{Summarizing tables}
\label{summarizing}

The tables of this section  summarize all fields in the database of degree
$\leq 11$.  
Numbers on tables which
are known to be correct are given in regular type.  Numbers which are
merely the bounds which come from perhaps incomplete lists of fields
are given in italics.
The table has a line for 
each group $nTj$, sorted by degree $n$ and the index $j$.    A more descriptive
name is given in the second column.

\begin{table}[htb]
\small
\setlength{\arraycolsep}{0.065in}
\[\begin{array}{|r|c|rrrr|rrr|r|}
\hline
\multicolumn{10}{|c|}{\textrm{Degree}\ 2}\\
\hline
T & G & \{2, 3\} & \{2, 5\} & \{3, 5\} & \{2, 3, 5\} &\textrm{rd}(K) & \textrm{grd}(K) & |\mathcal{K}[G,\Omega]| & \textrm{Tot} \\
\hline
1 & 2 & 7 & 7 & 3 & 15 & 1.73 & 1.73 & 1220 & 1216009 \\ 
\hline
\multicolumn{10}{c}{} \\[.08in]
\hline
\multicolumn{10}{|c|}{\textrm{Degree}\ 3}\\
\hline
T & G & \{2, 3\} & \{2, 5\} & \{3, 5\} & \{2, 3, 5\} &\textrm{rd}(K) & \textrm{grd}(K) & |\mathcal{K}[G,\Omega]| & \textrm{Tot} \\
\hline
1 & 3 & 1 & 0 & 1 & 1 & 3.66 & 3.66 & 47 & 1004 \\ 
\hline
2 & S_3 & 8 & 1 & 5 & 31 & 2.84 & 4.80 & 610 & 856373 \\ 
\hline
\multicolumn{10}{c}{} \\[.08in]
\hline
\multicolumn{10}{|c|}{\textrm{Degree}\ 4}\\
\hline
T & G & \{2, 3\} & \{2, 5\} & \{3, 5\} & \{2, 3, 5\} &\textrm{rd}(K) & \textrm{grd}(K) & |\mathcal{K}[G,\Omega]| & \textrm{Tot} \\
\hline
1 & 4 & 4 & 12 & 2 & 24 & 3.34 & 3.34 & 228 & 9950 \\ 
\hline
2 & 2^2 & 7 & 7 & 1 & 35 & 3.46 & 3.46 & 2421 & 52469 \\ 
\hline
3 & D_4 & 28 & 24 & 0 & 176 & 3.29 & 6.03 & 2850 & 764341 \\ 
\hline
4 & A_4 & 1 & 0 & 0 & 1 & 7.48 & 10.35 & 59 & 28786 \\ 
\hline
5 & S_4 & 22 & 3 & 1 & 143 & 3.89 & 13.56 & 527 & 281089 \\ 
\hline
\multicolumn{10}{c}{} \\[.08in]
\hline
\multicolumn{10}{|c|}{\textrm{Degree}\ 5}\\
\hline
T & G & \{2, 3\} & \{2, 5\} & \{3, 5\} & \{2, 3, 5\} &\textrm{rd}(K) & \textrm{grd}(K) & |\mathcal{K}[G,\Omega]| & \textrm{Tot} \\
\hline
1 & 5 & 0 & 1 & 1 & 1 & 6.81 & 6.81 & 7 & 181 \\ 
\hline
2 & D_5 & 0 & 4 & 2 & 8 & 4.66 & 6.86 & 146 & 11516 \\ 
\hline
3 & F_5 & 1 & 19 & 7 & 82 & 8.11 & 11.08 & 102 & 1645 \\ 
\hline
4 & A_5 & 0 & 5 & 6 & 62 & 7.14 & 18.70 & 78 & 95337 \\ 
\hline
5 & S_5 & 5 & 38 & 22 & 1353 & 4.38 & 24.18 & 192 & 598542 \\ 
\hline
\end{array}\]
\end{table}
\begin{table}[htb]
\small
\setlength{\arraycolsep}{0.05in}
\[\begin{array}{|r|c|rrrr|rrr|r|}
\hline
\multicolumn{10}{|c|}{\textrm{Degree}\ 6}\\
\hline
T & G & \{2, 3\} & \{2, 5\} & \{3, 5\} & \{2, 3, 5\} &\textrm{rd}(K) & \textrm{grd}(K) & |\mathcal{K}[G,\Omega]| & \textrm{Tot} \\
\hline
1 & 6 & 7 & 0 & 3 & 15 & 5.06 & 5.06 & 399 & 5291 \\ 
\hline
2 & S_3 & 8 & 1 & 5 & 31 & 4.80 & 4.80 & 610 & 8353 \\ 
\hline
3 & D_6 & 48 & 6 & 10 & 434 & 4.93 & 8.06 & 3590 & 147965 \\ 
\hline
4 & A_4 & 1 & 0 & 0 & 1 & 7.32 & 10.35 & 59 & 1357 \\ 
\hline
5 & 3 \wr 2 & 8 & 0 & 5 & 31 & 4.62 & 10.06 & 254 & 2169 \\ 
\hline
6 & 2 \wr 3 & 7 & 0 & 0 & 15 & 5.61 & 12.31 & 243 & 62484 \\ 
\hline
7 & S_4^+ & 22 & 3 & 1 & 143 & 5.69 & 13.56 & 527 & 242007 \\ 
\hline
8 & S_4 & 22 & 3 & 1 & 143 & 6.63 & 13.56 & 527 & 18738 \\ 
\hline
9 & S_3^2 & 22 & 0 & 4 & 375 & 7.89 & 15.53 & 445 & 9721 \\ 
\hline
10 & 3^2:4 & 4 & 0 & 2 & 44 & 8.98 & 23.57 & 34 & 396 \\ 
\hline
11 & 2 \wr S_3 & 132 & 18 & 2 & 2002 & 4.65 & 16.13 & 2196 & 323148 \\ 
\hline
12 & PSL_2(5) & 0 & 5 & 6 & 62 & 8.12 & 18.70 & 78 & 275 \\ 
\hline
13 & 3^2:D_4 & 50 & 0 & 0 & 624 & 4.76 & 21.76 & 274 & 27049 \\ 
\hline
14 & PGL_2(5) & 5 & 38 & 22 & 1353 & 11.01 & 24.18 & 192 & 11519 \\ 
\hline
15 & A_6 & 8 & 2 & 4 & \textit{540} & 8.12 & 31.66 & 10 & 670 \\ 
\hline
16 & S_6 & 54 & 30 & 42 & \textit{8334} & 4.95 & 33.50 & 26 & 21594 \\ 
\hline
\multicolumn{10}{c}{} \\[.08in]
\hline
\multicolumn{10}{|c|}{\textrm{Degree}\ 7}\\
\hline
T & G & \{2, 3\} & \{2, 5\} & \{3, 5\} & \{2, 3, 5\} &\textrm{rd}(K) & \textrm{grd}(K) & |\mathcal{K}[G,\Omega]| & \textrm{Tot} \\
\hline
1 & 7 & 0 & 0 & 0 & 0 & 17.93 & 17.93 & 4 & 117 \\ 
\hline
2 & D_7 & 0 & 0 & 0 & 0 & 6.21 & 8.43 & 80 & 496 \\ 
\hline
3 & 7:3 & 0 & 0 & 0 & 0 & 21.03 & 31.64 & 2 & 56 \\ 
\hline
4 & 7:6 & 0 & 0 & 1 & 5 & 12.10 & 15.99 & 94 & 189 \\ 
\hline
5 & SL_3(2) & 0 & 0 & 0 & \textit{} & 7.95 & 32.25 & \textit{36} & 618 \\ 
\hline
6 & A_7 & 0 & 2 & 3 & \textit{204} & 8.74 & 39.52 & \textit{1} & 331 \\ 
\hline
7 & S_7 & 10 & 24 & 14 & \textit{4391} & 5.65 & \textit{40.49} & \textit{1} & 8357 \\ 
\hline
\end{array}\]
\end{table}
\begin{table}[htbp]
\small
\setlength{\arraycolsep}{0.045in}
\[\begin{array}{|r|c|rrrr|rrr|r|}
\hline
\multicolumn{10}{|c|}{\textrm{Degree}\ 8}\\
\hline
T & G & \{2, 3\} & \{2, 5\} & \{3, 5\} & \{2, 3, 5\} &\textrm{rd}(K) & \textrm{grd}(K) & |\mathcal{K}[G,\Omega]| & \textrm{Tot} \\
\hline
1 & 8 & 4 & 8 & 0 & 16 & 11.93 & 11.93 & 23 & 5777 \\ 
\hline
2 & 4\times 2 & 6 & 18 & 1 & 84 & 5.79 & 5.79 & 581 & 15412 \\ 
\hline
3 & 2^3 & 1 & 1 & 0 & 15 & 6.93 & 6.93 & 908 & 10687 \\ 
\hline
4 & D_4 & 14 & 12 & 0 & 88 & 6.03 & 6.03 & 1425 & 24370 \\ 
\hline
5 & Q_8 & 2 & 0 & 0 & 8 & 18.24 & 18.24 & 7 & 778 \\ 
\hline
6 & D_8 & 20 & 20 & 0 & 104 & 6.71 & 9.75 & 708 & 29740 \\ 
\hline
7 & 8:\{1,5\} & 6 & 20 & 1 & 88 & 9.32 & 9.32 & 55 & 8040 \\ 
\hline
8 & 8:\{1,3\} & 22 & 10 & 0 & 120 & 10.09 & 10.46 & 121 & 10826 \\ 
\hline
9 & D_4\times 2 & 28 & 24 & 0 & 528 & 6.51 & 10.58 & 5908 & 175572 \\ 
\hline
10 & 2^2:4 & 8 & 24 & 0 & 160 & 6.09 & 9.46 & 620 & 29649 \\ 
\hline
11 & Q_8:2 & 18 & 18 & 0 & 312 & 6.51 & 9.80 & 921 & 17350 \\ 
\hline
12 & SL_2(3) & 0 & 0 & 0 & 0 & 12.77 & 29.84 & 4 & 681 \\ 
\hline
13 & A_4\times 2 & 7 & 0 & 0 & 15 & 8.06 & 12.31 & 243 & 26637 \\ 
\hline
14 & S_4 & 22 & 3 & 1 & 143 & 9.40 & 13.56 & 527 & 7203 \\ 
\hline
15 & 8:8^\times & 42 & 42 & 0 & 928 & 8.65 & 13.79 & 818 & 60490 \\ 
\hline
16 & 1/2[2^4]4 & 8 & 24 & 0 & 176 & 7.45 & 13.56 & 76 & 15545 \\ 
\hline
17 & 4\wr 2 & 16 & 72 & 0 & 480 & 5.79 & 13.37 & 252 & 42018 \\ 
\hline
18 & 2^2\wr 2 & 24 & 8 & 0 & 608 & 7.04 & 16.40 & 2544 & 216411 \\ 
\hline
19 & E(8):4 & 8 & 24 & 0 & 192 & 9.51 & 14.05 & 220 & 24380 \\ 
\hline
20 & [2^3]4 & 4 & 12 & 0 & 96 & 8.46 & 14.05 & 110 & 13631 \\ 
\hline
21 & 1/2[2^4]E(4) & 4 & 12 & 0 & 96 & 8.72 & 14.05 & 110 & 10091 \\ 
\hline
22 & E(8):D_4 & 0 & 0 & 0 & 204 & 8.43 & 18.42 & 882 & 19733 \\ 
\hline
23 & GL_2(3) & 128 & 24 & 4 & 912 & 8.31 & 16.52 & 388 & 6304 \\ 
\hline
24 & S_4 \times 2 & 132 & 18 & 2 & 2002 & 6.04 & 16.13 & 2196 & 45996 \\ 
\hline
25 & 2^3:7 & 0 & 0 & 0 & 0 & 12.50 & 17.93 & 1 & 20 \\ 
\hline
26 & 1/2[2^4]eD(4) & 64 & 24 & 0 & 1872 & 7.23 & 20.37 & 840 & 135840 \\ 
\hline
27 & 2\wr 4 & 16 & 48 & 0 & 448 & 5.95 & 19.44 & 160 & 86501 \\ 
\hline
28 & 1/2[2^4]dD(4) & 16 & 48 & 0 & 448 & 8.67 & 19.44 & 160 & 47150 \\ 
\hline
29 & E(8):D_8 & 48 & 24 & 0 & 1296 & 6.58 & 19.41 & 1374 & 170694 \\ 
\hline
30 & 1/2[2^4]cD(4) & 16 & 48 & 0 & 448 & 8.25 & 19.44 & 140 & 48249 \\ 
\hline
31 & 2\wr 2^2 & 16 & 8 & 0 & 432 & 5.92 & 19.41 & 458 & 54843 \\ 
\hline
32 & [2^3]A_4 & 0 & 0 & 0 & 0 & 13.56 & 34.97 & 24 & 29970 \\ 
\hline
33 & E(8):A_4 & 6 & 0 & 0 & 14 & 13.73 & 30.01 & 24 & 3240 \\ 
\hline
34 & E(4)^2:D_6 & 11 & 1 & 0 & 132 & 14.16 & 27.28 & 55 & 3907 \\ 
\hline
35 & 2\wr D(4) & 168 & 72 & 0 & 5568 & 5.83 & 22.91 & 1464 & 729730 \\ 
\hline
36 & 2^3:7:3 & 0 & 0 & 0 & 0 & 14.37 & 31.64 & 4 & 298 \\ 
\hline
37 & PSL_2(7) & 0 & 0 & \textit{} & \textit{} & \textit{21.00} & 32.25 & \textit{18} & 352 \\ 
\hline
38 & 2\wr A_4 & 24 & 0 & 0 & 112 & 10.66 & 37.27 & 46 & 67160 \\ 
\hline
39 & [2^3]S_4 & 168 & 24 & 0 & 2496 & 6.73 & 32.35 & 84 & 24625 \\ 
\hline
40 & 1/2[2^4]S(4) & 216 & 24 & 0 & 3872 & 7.67 & 29.71 & 98 & 12796 \\ 
\hline
41 & E(8):S_4 & 90 & 12 & 0 & 2282 & 8.38 & 28.11 & 222 & 11929 \\ 
\hline
42 & A_4\wr 2 & 12 & 0 & 0 & 83 & 7.68 & 32.18 & 14 & 3432 \\ 
\hline
43 & PGL_2(7) & \textit{4} & \textit{} & \textit{} & \textit{8} & \textit{11.96} & \textit{27.35} & \textit{27} & 1495 \\ 
\hline
44 & 2 \wr S_4 & 656 & 96 & 0 & 22944 & 5.84 & 31.38 & 336 & 440683 \\ 
\hline
45 & [1/2.S_4^2]2 & 110 & 0 & 0 & \textit{836} & 9.28 & 29.35 & 39 & 7732 \\ 
\hline
46 & 1/2[S(4)^2]2 & 28 & 0 & 0 & \textit{54} & 11.35 & 49.75 & 0 & 216 \\ 
\hline
47 & S_4\wr 2 & 542 & 0 & 0 & \textit{2185} & 6.74 & 35.05 & 15 & 28765 \\ 
\hline
48 & 2^3:SL_3(2) & 0 & 0 & \textit{} & \textit{} & \textit{11.36} & \textit{39.54} & \textit{6} & 495 \\ 
\hline
49 & A_8 & \textit{2} & \textit{4} & \textit{1} & \textit{55} & \textit{15.24} & \textit{72.03} & \textit{} & 90 \\ 
\hline
50 & S_8 & \textit{72} & \textit{30} & \textit{9} & \textit{1728} & \textit{11.33} & \textit{43.99} & \textit{1} & 4026 \\ 
\hline
\end{array}\]
\end{table}
\begin{table}[htbp]
\small
\setlength{\arraycolsep}{0.04in}
\[\begin{array}{|r|c|rrrr|rrr|r|}
\hline
\multicolumn{10}{|c|}{\textrm{Degree}\ 9}\\
\hline
T & G & \{2, 3\} & \{2, 5\} & \{3, 5\} & \{2, 3, 5\} &\textrm{rd}(K) & \textrm{grd}(K) & |\mathcal{K}[G,\Omega]| & \textrm{Tot} \\
\hline
1 & 9 & 1 & 0 & 1 & 1 & 13.70 & 13.70 & 3 & 52 \\ 
\hline
2 & 3^2 & 0 & 0 & 0 & 0 & 15.83 & 15.83 & 9 & 189 \\ 
\hline
3 & D_9 & 6 & 0 & 4 & 20 & 9.72 & 12.19 & 105 & 705 \\ 
\hline
4 & S_3 \times 3 & 8 & 0 & 5 & 31 & 8.38 & 10.06 & 254 & 10139 \\ 
\hline
5 & 3^2:2 & 1 & 0 & 1 & 15 & 14.29 & 15.19 & 48 & 373 \\ 
\hline
6 & 1/3[3^3]3 & 0 & 0 & 0 & 0 & 17.63 & 31.18 & 2 & 85 \\ 
\hline
7 & 3^2:3 & 0 & 0 & 0 & 0 & 26.09 & 50.20 & 0 & 90 \\ 
\hline
8 & S_3 \times S_3 & 22 & 0 & 4 & 375 & 8.93 & 15.53 & 445 & 7055 \\ 
\hline
9 & E(9):4 & 2 & 0 & 1 & 22 & 19.92 & 23.57 & 17 & 142 \\ 
\hline
10 & [3^2]S(3)_6 & 22 & 0 & 17 & 171 & 9.57 & 17.01 & 69 & 1066 \\ 
\hline
11 & E(9):6 & 6 & 0 & 4 & 20 & 14.67 & 16.83 & 64 & 880 \\ 
\hline
12 & [3^2]S(3) & 12 & 0 & 12 & 180 & 8.92 & 16.72 & 148 & 13929 \\ 
\hline
13 & E(9):D_6 & 6 & 0 & 4 & 20 & 10.98 & 16.83 & 64 & 642 \\ 
\hline
14 & 3^2:Q_8 & 4 & 0 & 0 & 19 & 21.52 & 29.72 & 2 & 47 \\ 
\hline
15 & E(9):8 & 5 & 1 & 0 & 18 & 17.74 & 25.41 & 3 & 40 \\ 
\hline
16 & E(9):D_8 & 25 & 0 & 0 & 312 & 9.19 & 21.76 & 137 & 434 \\ 
\hline
17 & 3\wr3 & 0 & 0 & 0 & 0 & 14.93 & 75.41 & 0 & 1274 \\ 
\hline
18 & E(9):D_{12} & 80 & 0 & 8 & 1380 & 8.53 & 22.06 & 290 & 9260 \\ 
\hline
19 & E(9):2D_8 & 60 & 1 & 0 & \textit{124} & 17.89 & 23.41 & 33 & 624 \\ 
\hline
20 & 3\wr S_3 & 18 & 0 & 12 & 60 & 7.83 & 29.89 & 30 & 7989 \\ 
\hline
21 & 1/2.[3^3:2]S_3 & 54 & 0 & 54 & 1296 & 9.82 & 24.90 & 126 & 4282 \\ 
\hline
22 & [3^3:2]3 & 18 & 0 & 12 & 60 & 10.27 & 26.46 & 51 & 784 \\ 
\hline
23 & E(9):2A_4 & 0 & 0 & 0 & 0 & 16.48 & 49.57 & 0 & 40 \\ 
\hline
24 & [3^3:2]S(3) & 321 & 0 & 48 & 8307 & 9.15 & 30.64 & 111 & 17973 \\ 
\hline
25 & [1/2.S(3)^3]3 & 4 & 0 & 0 & 4 & 12.89 & 29.96 & 4 & 303 \\ 
\hline
26 & E(9):2S_4 & 250 & 2 & 10 & \textit{362} & 12.79 & 27.88 & 51 & 866 \\ 
\hline
27 & PSL_2(8) & \textit{4} & \textit{} & \textit{} & \textit{4} & \textit{16.25} & \textit{30.31} & \textit{15} & 19 \\ 
\hline
28 & S_3 \wr3 & 28 & 0 & 0 & 90 & 8.18 & 33.56 & 7 & 6738 \\ 
\hline
29 & [1/2.S(3)^3]S(3) & 45 & 0 & 1 & 512 & 9.38 & 40.81 & 2 & 1255 \\ 
\hline
30 & 1/2[S(3)^3]S(3) & 232 & 1 & 40 & \textit{1637} & 6.86 & 30.37 & 35 & 5026 \\ 
\hline
31 & S_3\wr S_3 & 616 & 0 & 5 & 19865 & 6.83 & 36.26 & 15 & 112887 \\ 
\hline
32 & \Sigma L_2(8) & \textit{64} & \textit{} & \textit{} & \textit{240} & \textit{16.09} & \textit{34.36} & \textit{15} & 1141 \\ 
\hline
33 & A_9 & \textit{13} & \textit{} & \textit{2} & \textit{314} & \textit{14.17} & \textit{62.12} & \textit{} & 627 \\ 
\hline
34 & S_9 & \textit{46} & \textit{1} & \textit{1} & \textit{1507} & \textit{7.84} & \textit{53.19} & \textit{} & 3189 \\ 
\hline
\end{array}\]
\end{table}
\begin{table}[htb]
\small
\setlength{\arraycolsep}{0.045in}
\[\begin{array}{|r|c|rrrr|rrr|r|}
\hline
\multicolumn{10}{|c|}{\textrm{Degree}\ 10}\\
\hline
T & G & \{2, 3\} & \{2, 5\} & \{3, 5\} & \{2, 3, 5\} &\textrm{rd}(K) & \textrm{grd}(K) & |\mathcal{K}[G,\Omega]| & \textrm{Tot} \\
\hline
1 & 10 & 0 & 7 & 3 & 15 & 8.65 & 8.65 & 69 & 360 \\ 
\hline
2 & D_5 & 0 & 4 & 2 & 8 & 6.86 & 6.86 & 146 & 822 \\ 
\hline
3 & D_{10} & 0 & 24 & 4 & 112 & 8.08 & 10.91 & 768 & 857 \\ 
\hline
4 & 1/2[F(5)]2 & 1 & 19 & 7 & 82 & 10.23 & 11.08 & 102 & 178 \\ 
\hline
5 & F_5 \times 2 & 6 & 114 & 14 & 1148 & 9.48 & 14.50 & 584 & 1611 \\ 
\hline
6 & [5^2]2 & 0 & 8 & 4 & 16 & 6.84 & 18.02 & 32 & 175 \\ 
\hline
7 & A_5 & 0 & 5 & 6 & 62 & 12.35 & 18.70 & 78 & 146 \\ 
\hline
8 & [2^4]5 & 0 & 3 & 0 & 3 & 12.75 & 24.98 & 18 & 36 \\ 
\hline
9 & [1/2.D(5)^2]2 & 0 & 12 & 2 & 56 & \textit{12.71} & 24.72 & 34 & 87 \\ 
\hline
10 & 1/2[D(5)^2]2 & 0 & 22 & 12 & 126 & \textit{14.02} & 24.00 & 22 & 144 \\ 
\hline
11 & A_5 \times 2 & 0 & 35 & 18 & 930 & 9.42 & 22.24 & 179 & 1177 \\ 
\hline
12 & S_5(10a) & 5 & 38 & 22 & 1353 & 12.04 & 24.18 & 192 & 1560 \\ 
\hline
13 & S_5(10d) & 5 & 38 & 22 & 1353 & 9.16 & 24.18 & 192 & 1712 \\ 
\hline
14 & 2\wr 5 & 0 & 21 & 0 & 45 & 9.32 & 26.08 & 45 & 2050 \\ 
\hline
15 & [2^4]D(5) & 0 & 60 & 0 & 360 & 9.33 & 25.15 & 72 & 620 \\ 
\hline
16 & 1/2[2^5]D(5) & 0 & 60 & 0 & 360 & 9.46 & 25.15 & 72 & 509 \\ 
\hline
17 & [5^2:4]2 & 0 & 59 & 0 & 916 & \textit{17.46} & 26.65 & 34 & 1094 \\ 
\hline
18 & [5^2:4]2_2 & 0 & \textit{16} & 0 & \textit{17} & \textit{19.75} & 35.98 & 3 & 22 \\ 
\hline
19 & [5^2:4_2]2 & 0 & 17 & 0 & \textit{63} & \textit{16.96} & 28.08 & 18 & 111 \\ 
\hline
20 & [5^2:4_2]2_2 & 0 & 0 & 0 & \textit{31} & \textit{27.36} & 48.25 & 0 & 43 \\ 
\hline
21 & D_5\wr 2 & 0 & 34 & 0 & \textit{118} & 7.54 & 28.08 & 36 & 235 \\ 
\hline
22 & S_5\times 2 & 30 & 228 & 44 & 18942 & 7.06 & 26.99 & 570 & 26851 \\ 
\hline
23 & 2\wr D_5 & 0 & 360 & 0 & 5040 & 7.26 & 26.26 & 240 & 24024 \\ 
\hline
24 & [2^4]F(5) & 7 & 173 & 0 & 1250 & 14.13 & 27.62 & 30 & 1491 \\ 
\hline
25 & 1/2[2^5]F(5) & 7 & 173 & 0 & 1250 & 13.84 & 27.62 & 30 & 1491 \\ 
\hline
26 & PSL_2(9) & 4 & 1 & 2 & \textit{270} & 20.20 & 31.66 & 5 & 334 \\ 
\hline
27 & [1/2.F_5^2]2 & 0 & \textit{56} & 0 & \textit{652} & \textit{13.40} & 40.43 & 18 & 1052 \\ 
\hline
28 & 1/2[F_5^2]2 & \textit{1} & \textit{24} & \textit{2} & \textit{57} & \textit{15.16} & \textit{32.71} & \textit{8} & 59 \\ 
\hline
29 & 2\wr F_5 & 42 & 1038 & 0 & 17500 & 11.44 & 32.17 & 90 & 19112 \\ 
\hline
30 & PGL_2(9) & \textit{11} & \textit{5} & \textit{1} & \textit{55} & \textit{22.64} & \textit{34.42} & \textit{6} & 149 \\ 
\hline
31 & M_{10} & \textit{20} & \textit{4} & \textit{13} & \textit{83} & \textit{27.73} & \textit{53.50} & \textit{} & 198 \\ 
\hline
32 & S_6 & 27 & 15 & 21 & \textit{4166} & \textit{14.74} & 33.50 & 13 & 6913 \\ 
\hline
33 & F_5\wr 2 & 0 & \textit{177} & 0 & \textit{484} & 9.93 & \textit{35.41} & \textit{3} & 485 \\ 
\hline
34 & [2^4]A_5 & 0 & 35 & 0 & 1322 & 10.82 & 35.81 & 5 & 1388 \\ 
\hline
35 & P\Gamma L_2(9) & \textit{100} & \textit{32} & \textit{15} & \textit{1666} & \textit{17.98} & \textit{38.61} & \textit{15} & 3531 \\ 
\hline
36 & 2\wr A_5 & 0 & 245 & 0 & 19830 & 10.39 & 36.60 & 8 & 20660 \\ 
\hline
37 & [2^4]S_5 & 91 & 450 & 8 & 42059 & 7.80 & 38.11 & 17 & 60029 \\ 
\hline
38 & 1/2[2^5]S_5 & 91 & 450 & 8 & 42059 & 9.41 & 38.11 & 17 & 42851 \\ 
\hline
39 & 2\wr S_5 & 546 & 2700 & 16 & 588826 & 6.79 & 38.11 & 30 & 1095840 \\ 
\hline
40 & A_5\wr 2 & \textit{12} & \textit{63} & \textit{29} & \textit{1093} & 9.48 & \textit{41.90} & \textit{1} & 1124 \\ 
\hline
41 & [A_5:2]2 & \textit{28} & \textit{139} & \textit{11} & \textit{9435} & 9.30 & \textit{43.89} & \textit{1} & 9677 \\ 
\hline
42 & 1/2[S_5^2]2 & \textit{18} & \textit{68} & \textit{31} & \textit{850} & \textit{14.35} & \textit{45.93} & \textit{} & 882 \\ 
\hline
43 & S_5\wr 2 & \textit{185} & \textit{422} & \textit{} & \textit{20743} & 6.82 & \textit{48.97} & \textit{} & 31847 \\ 
\hline
44 & A_{10} & \textit{23} & \textit{16} & \textit{6} & \textit{801} & \textit{19.37} & \textit{51.68} & \textit{} & 1201 \\ 
\hline
45 & S_{10} & \textit{1} & \textit{12} & \textit{3} & \textit{2585} & \textit{7.77} & \textit{70.36} & \textit{} & 4944 \\ 
\hline
\end{array}\]
\end{table}
\begin{table}[htb]
\small
\setlength{\arraycolsep}{0.06in}
\[\begin{array}{|r|c|rrrr|rrr|r|}
\hline
\multicolumn{10}{|c|}{\textrm{Degree}\ 11}\\
\hline
T & G & \{2, 3\} & \{2, 5\} & \{3, 5\} & \{2, 3, 5\} &\textrm{rd}(K) & \textrm{grd}(K) & |\mathcal{K}[G,\Omega]| & \textrm{Tot} \\
\hline
1 & 11 & 0 & 0 & 0 & 0 & 17.30 & 17.30 & 1 & 18 \\ 
\hline
2 & D_{11} & 0 & 0 & 0 & 0 & 10.24 & 12.92 & 32 & 55 \\ 
\hline
3 & 11:5 & 0 & 0 & 0 & 0 & \textit{88.82} & \textit{105.74} & 0 & 2 \\ 
\hline
4 & 11:10 & \textit{} & \textit{} & \textit{} & \textit{} & \textit{17.01} & 20.70 & \textit{4} & 55 \\ 
\hline
5 & PSL_2(11) & \textit{} & \textit{} & \textit{} & \textit{} & \textit{15.36} & \textit{42.79} & \textit{2} & 91 \\ 
\hline
6 & M_{11} & \textit{} & \textit{} & \textit{} & \textit{1} & \textit{96.24} & \textit{270.83} & \textit{} & 10 \\ 
\hline
7 & A_{11} & \textit{} & \textit{} & \textit{} & \textit{4} & \textit{21.15} & \textit{146.24} & \textit{} & 71 \\ 
\hline
8 & S_{11} & \textit{5} & \textit{4} & \textit{1} & \textit{123} & \textit{7.72} & \textit{91.50} & \textit{} & 931 \\ 
\hline
\end{array}\]
\end{table}

The next four columns represent a main focus of the database, complete
lists of fields ramified within a given set of primes.  As a matter of
notation, we write e.g.\ $\cK(G,-^*p^*q^*)$ to denote the union of all
$\cK(G,-^sp^aq^b)$.  The database contains completeness results for
many other prime combinations beyond those given in the table;
\S\ref{quintics2357}-\S\ref{octicwild} give examples of these further
completeness results.

The next column gives minimal values of root discriminants.  
More refined minima can easily be obtained from the database.  For example, 
for $S_5$, minimal discriminants for $s = 0$, $1$, and $2$ complex places
are respectively $(61 \cdot 131)^{1/5} \approx 7.53$,  $(13 \cdot
347)^{1/5} \approx 5.38$, and  
$1609^{1/5} \approx 4.38$.  Completeness is typically known well
past the minimum.   

In understanding root discriminants, the Serre-Odlyzko constant 
$\Omega = 8 \pi e^\gamma \approx 44.76$  
mentioned in the introduction plays an important role as follows.   
First, if $K$ has root discriminant $<\Omega/2$, then its
maximal unramified  extension $K'$ has finite
degree over $K$.  Second, if $\rd(K) < \Omega$, then
the generalized Riemann hypothesis implies the same conclusion
$[K':K]<\infty$.  Third, suggesting that there is 
a modestly sharp qualitative transition associated with
$\Omega$, the field $\Q(e^{2 \pi i/81})$ with root 
discriminant  $3^{3.5} \approx 46.77$ has $[K':K] = \infty$ by
 \cite{hoelscher}.  

The next two columns of the tables again represent a main focus of the database,
complete lists of fields with small Galois root discriminant.   We
write $\cK [G,B]$ for the set of all fields with Galois group $G$ 
and  $\mbox{grd}(K) \leq B$.  The tables
give first the minimal Galois root discriminant.   They next give
$|\cK[ G,\Omega]|$.  
For many groups, the database is complete for cutoffs well past
$\Omega$.  For example, the set $\cK[9T17,\Omega]$ is empty, and not
adequate for the purposes of \cite{jr-Artin}.  However the database
identifies $|\cK[9T17,200]| = 36$ and this result is adequate for the
application.
 
The last column in a table gives the total number of fields in the database for
the given group.  Note that one could easily make this number much
larger in any case.  For example, a regular family over $\Q(t)$ for
each group is given in \cite[App.~1]{malle-matzat}, and one could
simply specialize at many rational numbers $t$.  However we do not
do this: all the fields on our database are there
only because discriminants met one criterion or another for being
small.  The fluctuations in this column should not be viewed as
significant, as the criteria depend on the group in ways driven
erratically by applications.
 
There are a number of patterns on the summarizing tables which hold
because of relations between transitive groups.  For example the
groups $5T4 = A_5$, $6T12 = PSL_2(5)$, and $10T7$ are all isomorphic.
Most of the corresponding lines necessarily agree.  Similarly $A_5$ is
a quotient of $10T11$, $10T34$, and $10T36$.  Thus the fact that
$\cK(A_5,-^*2^* 3^*) = \emptyset$ immediately implies that also
$\cK(10Tj,-^*2^* 3^*) = \emptyset$ for $j \in \{11,34,36\}$.

Almost all fields in the database come from complete searches of
number fields carried out by the authors.  In a few cases, we obtained
polynomials from other sources, notably for number fields of small
discriminant: those compiled by the Bordeaux group \cite{megrez},
which in turn were computed by several authors, and the the tables of
totally real fields of Voight \cite{voight-tr, voight-web}.  In
addition, we include fields found by the authors in joint work with
others \cite{driver-jones2, jj-wallington}.

To compute cubic fields, we used the program of Belabas \cite{belabas,
  belabas-prog}.  Otherwise, we obtain complete lists by using 
   traditional and targeted Hunter
searches \cite{jr1,jr-septic} or the class field theory functions in pari/gp
\cite{gp}.  For larger nonsolvable groups where completeness results
are currently out of reach,
we obtain most of our fields by specializations of families at carefully chosen
points to keep ramification small in various senses.

\section{$S_5$ quintics with discriminant $-^s 2^a 3^b 5^c 7^d$}
\label{quintics2357}
    One of our longest searches determined $\cK(S_5,-^*2^* 3^* 5^* 7^*)$, finding
it to consist of 11279 fields.  In this section, we consider how this set 
interacts with mass heuristics.  

In general, mass heuristics \cite{Bh,malle-dist} give one
expectations as to the sizes $|\cK(G,D)|$  
of the sets contained in the database.   Here we consider these
heuristics only in  
the most studied case $G=S_n$.    The mass of a $\Q_v$-algebra $K_v$
is by definition $1/|\mbox{Aut}(K_v)|$.  Thus the mass of $\R^{n-2s} \C^s$ is
\begin{equation}
\label{massinf}
\mu_{n,-^s} = \frac{1}{(n-2s)! s! 2^s}.
\end{equation}
For $p$ a prime, similarly let $\mu_{n,p^c}$ be the total
mass of all $p$-adic algebras with degree $n$ and discriminant
$p^c$.   For $n<p$, all algebras involved are tame and
\begin{equation}
\label{masstame}
\mu_{n,p^c} = |\{\mbox{Partitions of $n$ having $n-c$ parts}\}|.
\end{equation}
For $n \geq p$, wild algebras are involved.  General formulas for $\mu_{n,p^c}$ are given in 
 \cite{Ro-WP}.

The mass heuristic
says that if the discriminant $D=-^s \prod_p p^{c_p}$ in question is
a non-square, then 
\begin{equation}
\label{heur}
|\cK(S_n,-^s \prod_p p^{c_p})| \approx  \delta_n \mu_{n,-^s} \prod_p \mu_{n,p^{c_p}}.
\end{equation}
Here $\delta_n = 1/2$, except for the special cases $\delta_1=\delta_2=1$ which
require adjustment for simple reasons.  
The left side is an integer and the right side is often 
close to zero because of \eqref{massinf} and \eqref{masstame}.
So \eqref{heur} is intended only to be 
used in suitable averages.

For $n \leq 5$ fixed and $|D| \rightarrow \infty$,
the  heuristics are exactly right on average, the case $n=3$ being the
Davenport-Heilbronn theorem and the cases $n=4$, $5$ more
recent results of Bhargava \cite{bhargava4, bhargava5}.
For a fixed set of ramifying primes $S$ and $n \rightarrow \infty$,
the mass heuristic predicts no fields after a fairly sharp 
cutoff $N(S)$, while in fact there can be many fields in degrees 
well past this cutoff \cite{RobCheb}.   Thus the 
regime of applicability of the mass heuristic is 
not clear.  

To get a better understanding of this regime,
it is of interest to consider other limits.   Let $c_n$ be the number
of elements of order $\leq 2$ in $S_n$.  Let $\mu_{n,p^*}$ be the 
total mass of all $\Q_p$-algebras of degree $n$.  Thus $\mu_{n,p^*}$ is
the number of partitions of $n$ if $n<p$.   Then, for $k \rightarrow \infty$, 
the mass heuristic predicts the asymptotic equivalence
\begin{equation}
\label{massprimecut}
|\cK(S_n,-^* 2^* \cdots p_k^*\})| \sim \delta_n \mu_{n,-^*} \prod_{j=1}^k \mu_{n,p^*}.
\end{equation}
Both sides of \eqref{massprimecut} are $1$ for all $k$ when $n=1$.  
For $n=2$ and $k \geq 1$, the statement becomes $2^{k+1}-1 \sim 2^{k+1}$ which is true.  
Using the fields in the database as a starting point, we have carried out substantial
calculations suggesting that, after removing fields with discriminants of the 
form $-3u^2$ from the count on the left, \eqref{massprimecut} holds also 
for $n=3$ and $n=4$.  

In this section, we focus on the first nonsolvable case, $n=5$.   For $k \geq 3$, \eqref{massprimecut}
becomes 
\begin{equation}
\label{mass5}
|\cK(S_5,-^* 2^* \cdots p_k^*)|  \stackrel{}{\sim} \frac{1}{2} \cdot \frac{26}{120} \cdot 40 \cdot 19 \cdot 27 \cdot 7^{k-3} \approx 6.48 \cdot 7^k.
\end{equation}
Through the cutoff $k=4$, there are fewer fields than predicted by the
mass heuristic:
\[
\begin{array}{r| cccc}
p_k    &                    2 & 3 & 5 & 7  \\
\hline
|\cK(S_5,-^* 2^* \cdots p_k^*)|  &                    0 & 5 & 1353 & 11279 \\
                                                           &                                0\% & 6\% & 61\% & 72\% 
 \end{array}.
\]
For comparison, the ratio $11279/(6.48 \cdot 7^4) \approx 72\%$ is actually larger than the ratios at $k=4$ 
for cubic and quartic fields with discriminants $-3u^2$ removed, 
these being respectively $64/(1.33 \cdot 3^4) \approx 47\%$ and $740/(3.30 \cdot 5^k) \approx 36\%$.  
As remarked above, these other cases experimentally approach $100\%$ as $k$ increases.  This
experimental finding lets one
reasonably argue that \eqref{mass5} may hold too, with the small percentage
72\% being a consequence of a small discriminant effect.

\begin{table}[htb]
{\renewcommand{\arraycolsep}{3.5pt}
\[
\begin{array}{c|llllllllllll|r}
v \setminus c & 0 & 1 & 2 & 3 & 4 & 5 & 6 & 7 & 8 & 9 & 10 & 11 & \mbox{Total} \\
\hline
\infty & 1 & 10 & 15 & &&&&&&&&&26\\
         &0.71 & 9.52 & 15.77 &  &&&&&&&&&    \\
         \hline
2 & 1 &   & 2 & 2 & 5 & 4 & 6 &  & 4 & 4 & 4 & 8 & 40 \\
 & 0.73 &  & 1.66 & 1.48 & 4.71 & 3.83 & 5.66 &  & 4.47 & 4.37 & 4.15 & 8.94 & \\
\hline
3 & 1 & 1 & 1 &  3 & 5 & 5 & 3 &&&&&&19 \\
  & 0.76 & 0.85 & 0.78 & 2.89 & 5.24 & 5.13 & 3.43 &&&&&\\
\hline
5 & 1 & 1 & 2 & 2 &  & 4 & 4 & 4 & 4 & 5&&& 27 \\
   & 0.37 & 0.39 & 0.96 & 1.32 &    & 4.07 & 4.17 & 4.70 & 4.65 & 6.38 & & & \\
\hline
7 & 1 & 1 & 2 & 2 & 1 & &&&&&&&7 \\
   & 0.84& 0.88 &1.92 & 2.12 & 1.24 & \\
\end{array}
\]
}
\caption{\label{quintrel}  Local masses $120 \mu_{5,-^c}$ and 
$\mu_{5,p^c}$, compared with frequencies of local
discriminants from $\cK(S_5,-^*2^*3^*5^*7^*)$.  }
\end{table}

Table~\ref{quintrel} compares local masses 
with frequencies of actually occurring local discriminants, 
inflated by the ratio $(6.58 \cdot 7^4)/11279$ to facilitate direct 
comparison.   Thus, e.g., the $7$-adic discriminants $(7^0, 7^1, 7^2, 7^3, 7^4)$ are
predicted by the mass heuristic to occur with relatively frequency
$(1,1,2,2,1)$.  They
actually occur with relative frequency $(0.84, 0.88, 1.92, 2.12, 1.24)$.   
Here and for the other four places, trends away from the predicted
asymptotic values are explained by consistent underrepresentation
of fields with small discriminant.  The consistency of the data with the
mass heuristic on this refined level provides further
support for \eqref{mass5}.

\section{Low degree nonsolvable fields with discriminant $-^s p^a q^b$}
\label{nonsolvable567}
Our earliest contributions to the general subject of number field
tabulation were \cite{jr1} and \cite{jr-septic}. These papers respectively 
found that there are exactly $398$ sextic and $10$ septic fields
with discriminant of the form $-^s 2^a 3^b$.   On the lists from these papers, the 
nonsolvable groups $PSL_2(5)  \cong A_5$, $PGL_2(5) \cong A_6$, $A_6$, $S_6$ and $SL_3(2)$, $A_7$, $S_7$
respectively arise $0$, $5$, $8$, $54$, and $0$, $0$, $10$ times.  In this section, 
we summarize further results from the database of this form, identifying
or providing lower bounds for $|\cK(G, -^* p^*q^*)|$.    

\begin{table}[htb]
{\renewcommand{\arraycolsep}{1.5pt}
\[
\begin{array}{r|rrr|rrrrrrrrrrrrrrrrrrrrrr|r}
 &\; 2&\; 3&\; 5& \; 7&11&13&17&19&23&29&31&37&41&43&47&53&59&61&67&71&73&79&83&89&97 & T \\
 \hline
2& \bu    &5&38&2&2&4&3 &2&5&6&3&6 &9 &14&10&11&8&13&13&8 &5&11&10&13&4 & 205  \\
3&     & \bu &22& & & &1 &4& & &1&2 &2 &6 &3 &5 &3&3 &2 &8 &3&2 &4 &3 &2 & 81  \\
5&    5&6& \bu  &4&5&9&12&8&8&8&9&13&12&11&8 &14&8&15&13&14&9&14&11&11&14 & 290 \\
    \hline
7&     & &  & \bu&1& &  & &1&2& &  &2 &2 &  &  & &  &1 &3 & &1 &  &1 & & 20   \\
11&     &2&  & &\bu & &  &1&1& & &1 &  &1 &1 &  &1&  &  &1 &1&  &1 &1 &  & 18\\
13&     & &1 & & &\bu &  &1& &2& &3 &2 &  &1 &  & &  &  &  &1&  &  &1 &1 & 25 \\
17&    1&1&  & & & &\bu  &1& &3& &3 &4 &1 &  &1 &1&  &  &1 &2&  &  &1 &3 & 37\\
 19&   1&3&  & & & &  & \bu &2& &3&  &1 &  &1 &1 &1&1 &3 &2 & &1 &1 &  &2 & 36  \\
23&     &1&  & & & &  & & \bu&1& &1 &  &2 &  &  &1&1 &  &  & &  &2 &1 &2 & 28 \\
29&     2&3&  & & &1&2 & &2& \bu& &2 &3 &1 &1 &  &2&3 &3 &1 &4&2 &  &2 &2 & 48  \\
31&    1&3&1 & & &2&1 & & & & \bu&1 &  &4 &1 &  & &2 &2 &1 & &  &2 &  &1 & 30  \\
37&     & &1 & & & &  & & & & &  \bu&2 &3 &  &2 &2&2 &  &  & &3 &4 &1 &1 & 52 \\
41&   2&2&2 & & & &1 &1&1& & &  & \bu  &1 &  &2 & &3 &  &1 &1&2 &1 &3 &5 & 56 \\
43&    1&3&1 & &1& &  & & & & &  &1 & \bu &  &  &2&1 &1 &2 &2&4 &1 &2 &  & 61 \\
47&     & &  & & & &  & & & & &  &  &  & \bu &1 &4&  &3 &  & &1 &3 &  & & 38  \\
53&     & &2 & &1& &  & &1& & &  &  &  &  & \bu &3&3 &  &1 &2&1 &3 &2 & & 52 \\
59&   1&3&2 & & &1&2 & & & &1&  &  &  &  &  & \bu&  &  &1 &1&2 &1 &3 &1&45 \\
 61&   1&1&  & & & &1 & & &1& &1 &1 &  &  &  & &  \bu &2 &3 &2&1 &  &  &1 & 56\\
 67&   2&1&1 & & & &  & & & & &  &  &  &1 &  & &  &  \bu &4 & &2 &2 &1 &2 & 54\\
 71 &   1&2&1 & & & &  &1& &1&4&2 &  &  &  &  &2&1 &  & \bu &  &4 &1 &3 & &59 \\
 73 &   1& &1 &1& & &  & & &1& &  &2 &1 &1 &  & &1 &  &  &\bu &  &  &2 &3 & 38 \\
 79&    4&4&2 & & &1&  &2& & & &  &3 &  &  &1 &1&  &  &  & &\bu  &  &2 &1 & 54 \\
 83&    &1&  & & & &  & & & & &  &  &  &  &1 & &  &1 &  & &2 & \bu &4 &  & 51 \\
 89&   1&3&2 & & & &  & & &1&1&1 &2 &  &  &2 & &1 &1 &1 & &2 &  &\bu  &1 & 58 \\
 97&    &1&  & & &1&  & & & & &  &1 &  &  &  &1&  &1 &  & &1 &  &1 &  \bu & 46 \\
 \hline
 T & 24 & 40 & 28 & 1 & 4 & 7 & 9 & 8 & 5 & 14 & 14 & 5 & 19 & 8 & 2 & 8 & 14 & 9 & 8 & 16 & 9 & 23 & 5 & 19 & 7 \\
     \end{array}
\]
}
\caption{\label{quintictable} $|\cK(A_5,-^* p^*q^*)|$ beneath the diagonal and $|\cK(S_5,-^* p^* q^*)|$ above the diagonal.  
It is expected that $A_5$ totals are smaller for primes $p \equiv 2,3 \; (5)$ because in this case 
$p^4$ is not a possible local discriminant. }
\end{table}

The format of our tables exploits the fact that in the range
considered for a given group, there are no fields ramified at one prime only. 
In fact \cite{jr-oneprime}, the smallest prime $p$ for which $\cK(G,-^* p^*)$ is nonempty is
as follows:
\[
\begin{array}{c||cc|cc|ccc|c}
G &  A_5 & S_5 & A_6 & S_6 & SL_3(2) & A_7 & S_7 & PGL_2(7) \\
\hline
p   & 653 & 101 & 1579 & 197 & 227 & \!\!\!> \!227 & 191 & 53 
\end{array}.
\]
Here $7T5 = SL_3(2)$ is abstractly isomorphic to $8T37 = PSL_2(7)$ and
thus has index two in $8T43 = PGL_2(7)$.

\begin{table}[htb]
\[
\begin{array}{r|rrr|rrrrrrrr|r}
&2 & 3 & 5 & 7 & 11 & 13 & 17 & 19 & 23 & 29 & 31 & T \\
\hline
2& \bu &54&30& & &2&2&2&4&4 &6 & 104\\
3&                     8& \bu  &42& &4& &8& &2&12&2 & 124 \\
5&                     2&4 & \bu &2&2&2& &6&8&2 &4 & 98\\
\hline
7&                      &2 &  &\bu & & & & & &  & & 2 \\
11&                     &  &  & &\bu & & & & &  & & 6  \\
13&                      &2 &  & & & \bu& & & &  & & 4 \\
17 &                      &  &2 & & & &\bu & &2&  & & 12  \\
19&                      &2 &  & & &2& &\bu & &  & & 8  \\
23&                     2&2 &  & & & & & & \bu &  & & 16 \\
29&                      &4 &  & & & & &2& & \bu & & 18 \\
31&                     4&6 &  & & & & & & &2 & \bu & 12 \\
\hline
T & 16 & 30 & 8 & 2 &  & 4 & 2 & 6 & 4 & 8 & 12 \\
\end{array}
\]
\caption{\label{sextictable} $|\cK(A_6,-^* p^*q^*)|$ beneath the diagonal and $|\cK(S_6,-^* p^*q^*)|$ above the diagonal.
All entries are even because contributing fields come in twin pairs. }
\end{table}

\begin{table}[htb]
\[
\begin{array}{r|rrrrrr}
\multicolumn{7}{c}{SL_3(2) \mbox{ and } PGL_2(7)} \\
    & 2 & 3 & 5 & 7 & 11 & 13 \\
    \hline
2& \bu    & \mathit{4 }    &   \mathit{0 }     &  \mathit{51 }  &  \mathit{0 }   &  \mathit{0 }   \\
3&   0                   & \bu    &    \mathit{ 0}   &  \mathit{28 }  &  \mathit{0 }   &  \mathit{ 0 }  \\
5&   0                  & 0& \bu    &    \mathit{4 } &  \mathit{0 } &  \mathit{ 0 }  \\
7&   \mathit{44}              & \mathit{12}  & \mathit{4}   & \bu  &  \mathit{ 4}  &  \mathit{ 6 }  \\
11&   4                    &  0 & \mathit{0}  & \mathit{6}  & \bu &  \mathit{ 0 }   \\
13&   0                     & \mathit{0}   & \mathit{0}   & \mathit{0} & 0 & \bu \\
\end{array}
\;\;\;\;\;\;\;\;\;\;\;\;\;
\begin{array}{r|rrrrrr}
\multicolumn{7}{c}{A_7 \mbox{ and } S_7} \\
    & 2 & 3 & 5 & 7 & 11 & 13 \\
    \hline
2& \bu    &  10    &   24  &  \mathit{55} & 0  & 0  \\
3&                   0  & \bu    &  14   & \mathit{44} & 2  &  \mathit{0} \\
5&                   2 & 3 & \bu    &  \mathit{18}  & \mathit{0} & \mathit{0}  \\
7&                \mathit{0}  & \mathit{7} & \mathit{5} & \bu  & \mathit{5}  & \mathit{0}\\
11&                    0    & 0    & \mathit{1}  & \mathit{0}  & \bu & 0  \\
13&                    0    & \mathit{0}    & \mathit{0}    & \mathit{0}  & 0 & \bu \\
\end{array}
\]
\caption{\label{septictable}Determinations or lower bounds for 
$|\cK(G,-^* p^*q^*)|$ for four $G$.   The entries 
$|\cK(SL_3(2),-^* p^*q^*)|$ are all even because contributing fields
come in twin pairs.    }
\end{table}

Restricting to the six groups $G$ of the form $A_n$ or $S_n$,  our
results on $|\cK(G,-^* p^*q^*)|$ compare with the mass heuristic as
follows.  First, local masses $\mu_{n,v^*}$ are given in the middle six columns:
 \begin{equation}
 \label{localmasses}
 \begin{array}{r | r @{/} l  r  r  r  r  r  | c } \\
 n & \multicolumn{2}{c}{ \infty } & 2 & 3 & 5 & 7 & \mbox{tame} & \mu \\
 \hline
 5 & 26 & 120 & 40 & 19 & 27 &  & 7 & 5.31  \\
 6 & 76 & 720 & 145 & 83 & 31 &  &  11 & 6.39 \\
 7 & 232 & 5040 & 180 & 99 & 55 & 57 & 15 & 5.18 
 \end{array}.
 \end{equation}
The column $\mu$ contains the global mass $0.5 \mu_{n,-^*} \mu_{n,p^*} \mu_{n,q^*}$ 
for two tame primes $p$ and $q$.   When one or both of the primes are wild,
the corresponding global mass is substantially larger.

Tables~\ref{quintictable}, \ref{sextictable}, and \ref{septictable}
clearly show that there tend to be more fields when one or more 
of the primes $p$, $q$ allows wild ramification, as one would expect from
\eqref{localmasses}.  To make plausible conjectures about the asymptotic
behavior of the numbers $|\cK(G,-^* p^* q^*)|$, one would have to do more complicated local
calculations  than those summarized in  \eqref{localmasses}.  These calculations would have
to take into account various secondary phenomena, such as the fact that $s$ is forced
to be even if $p \equiv q \equiv 1 \; (4)$.
Tables~\ref{quintictable}, \ref{sextictable}, and \ref{septictable}
each 
reflect substantial computation, but the amount of evidence is too
small to warrant making formal conjectures
in this setting.

\section{Nilpotent octic fields with odd discriminant $-^s p^a q^b$}
\label{octictame}  
   The database has all 
    octic fields with Galois group 
a $2$-group and discriminant of the form $-^s p^a q^b$ 
with $p$ and $q$ odd primes $<250$.  There are $\binom{52}{2} = 1326$
pairs $\{p,q\}$ and the average size of $\cK(\mbox{NilOct},-^*p^*q^*)$ in this
range is about $12.01$.   In comparison with the
nonsolvable cases discussed in the previous two sections, there
is much greater regularity in this setting.  We
exhibit some of the greater regularity and explain how it
makes  some of the abstract considerations
of \cite{boston-ellenberg,boston-perry} more concrete.

\begin{table}[htb]
{\small
{\renewcommand{\arraycolsep}{2.5pt}
\[
\begin{array}{|rr|rrrr|rrrr|rrrrr|rrr|r|r|rrr|}
\hline
p & q & 1 & 2 & 4 & 5 & 6 & 7 & 8 & 10 & 16 & 17 & 19 & 20 & 21 & 27 & 28 & 30 & s_3 & \# & \nu & T &  s \\
\hline
\multicolumn{20}{c}{\;} \\
\cline{1-2} \cline{21-21} 
                     3,7 & 3,7 & \multicolumn{18}{c}{\;} & \multicolumn{1}{|c|}{1/4}  \\
\hline
                     7_2 &3_1  &  & &1& & & &1& & &  &  & & &  &  &  &4 &193 & 1/8 & \circ & 4\\
                     11_2 &7_1 &  & &1& &2& & & & &  &  & & &  &  &  &4 &185 & 1/8 & \bullet & \geq 5 \\
\hline
\multicolumn{20}{c}{\;} \\
\cline{1-2} \cline{21-21} 
                     3,7 & 5 & \multicolumn{18}{c}{\;} & \multicolumn{1}{|c|}{1/4}\\
  \hline
 3_1 &5_1  &  &1& & & &1& & & &  &  & & &  &  &  &4 &219 & 1/8 & i & 4 \\
 11_4 &5_2 &  &1&2& & & & &2&2&4 &  & & &  &  &  &6 &87 & 1/16 & ii & 6\\
 19_2& 5_2 &  &1&2& &2& &1&2& &2 &2 &1&1&2 &2 &4 &7 &86 & 1/16 & iii & 19\\
\hline
\multicolumn{20}{c}{\;} \\
\cline{1-2} \cline{21-21} 
                     3,7 & 1 & \multicolumn{18}{c}{\;} & \multicolumn{1}{|c|}{1/4} \\
                     \hline
                       3_1&17_1 &2 &1& & & & & & & &  &  & & &  &  &  &4 &162 & 1/8 & & \\
                       19_2&17_2 &2 &1&2& &2&1&1&2& &2 &  & & &  &  &  &6 &66 & 1/16 & &  \\
                       23_4&41_2 &2 &1&2& &4&1&2&2&2&4 &2 &1&1&4 &4 &4 &8 &52 & 1/16 & &  \\
  \hline
\multicolumn{20}{c}{\;} \\
\cline{1-2} \cline{21-21} 
                     5 & 5 & \multicolumn{18}{c}{\;} & \multicolumn{1}{|c|}{1/16} \\
\hline
                     5_1 &13_1 &  &3& & & &2& & & &  &  & & &  &  &  &5 &42 & 1/32 & I  &  6\\
                     5_2 &29_2 &  &3&3&1& & & &6& &8 &4 &2&2&4 &4 &  &9 &11 & 1/128 & II & 27 \\
                     13_4&53_4 &  &3&3&1&2& &2&6& &8 &12&6&6&12&12&16&11&13 & 1/128 & III & \geq \! 17 \\
                     13_2&29_4 &  &3&3&1& & & &6&2&12&4 &2&2&2 &2 &4 &9 &25 & 1/64& IV & \geq \! 30 \\
\hline
\multicolumn{20}{c}{\;} \\
\cline{1-2} \cline{21-21} 
                     5 & 1 & \multicolumn{18}{c}{\;}& \multicolumn{1}{|c|}{1/8} \\
\hline
                  5_1 &17_1 &4 &3& & & & & & & &  &  & & &  &  &  &5 &76 & 1/16 & & \\
                    13_4&17_2 &4 &3&3&1& &2& &6&2&8 &4 &2&2&4 &4 &  &9 &17 & 1/64 & & \\
                  5_2 &41_4 &4 &3&3&1& &2& &6& &4 &4 &2&2&2 &2 &4 &9 &22 & 1/64 & &  \\
                  53_2 &17_2 &4 &3&3&1&2&2&2&6& &12&4 &2&2&2 &2 &4 &9 &18 & 1/64 & &  \\
                  109_4 & 73_4 &4 &3&3&1&4&2&4&6&2&16&12&6&6&12&12&16&11&6 & 1/128 & &  \\
                 101_4& 97_4&4 &3&3&1&4&2&4&6&2&16&12&6&6&16&16&24&12&1& 1/128 & &   \\

\hline
\multicolumn{20}{c}{\;} \\
\cline{1-2} \cline{21-21} 
                     1 & 1 & \multicolumn{18}{c}{\;} & \multicolumn{1}{|c|}{1/16}\\
\hline
                     17_1&41_1 &12&3& & & & & & & &  &  & & &  &  &  &6 &27 & 1/32 &  & \\
                     41_4&73_2 &12&3&3&1& &6& &6&2&  &4 &2&2&4 &4 &  &9 &12 & 1/64 & & \\
                     41_2 & 241_2 &  12 & 3 & 3 & 1 & 4 & 6 & 4 & 6 & 2 & 16 & 4 & 2 & 2 & 4 & 4 & 4 & 10 & 2 &  1/128 & &  \\
                     73_4&89_4 &12&3&3&1&6&6&6&6&6&24&12&6&6&16&16&16&12&2 & 1/256 & &   \\
                     73_4&137_4 &12&3&3&1&6&6&6&6&6&24&12&6&6&24&24&24&13&2 & 1/256 & & \\
\hline
\end{array}
\]
}
}
\caption{\label{pqoctictab} Nonzero cardinalities $|\cK(8Tj,-^* p^* q^*)|$ for $8Tj$ an octic group 
of $2$-power order }
\end{table}

Twenty-six of the fifty octic groups have $2$-power order. 
 Table~\ref{pqoctictab} presents the nonzero cardinalities, so that
e.g. $|\cK(8Tj,-^* 5^* 29^*)| = 4, 2, 2$ for $j=19$, $20$, $21$.  
The repeated proportion $(2,1,1)$ for these groups and 
other similar patterns are due to the sibling phenomenon discussed
in Section~\ref{summarizing}.    Only the sixteen $2$-groups
generated by $\leq 2$ elements actually occur.  Columns  $s_3$, $\#$, $\nu$, $T$,
and $s$ are all explained later in this section.

The main phenomenon presented in Table~\ref{pqoctictab}  is that 
the multiplicities presented are highly repetitious, with e.g.\ the 
multiplicities presented for $(5,29)$ occurring for all together
eleven pairs $(p,q)$, as indicated in the $\#$ column.   
The repetition is even greater than indicated by 
the table itself.  Namely if $(p_1,q_1)$ and $(p_2,q_2)$ 
correspond to the same line, then not only are
the numbers $\cK(8Tj,-^* p_i^* q_i^*)$ independent
of $i$, but the individual $\cK(8Tj,-^* p_i^a q_i^b)$ 
and even further refinements are also independent of $i$.  

The line corresponding to a given pair $(p,q)$ is almost determined by elementary considerations as follows.  Let
$U$ be the order of $q$ in $(\Z/p)^\times$ and let 
$V$ be the order of  $p$ in $(\Z/q)^\times$.  Let 
$u = \gcd(U,4)$ and $v = \gcd(V,4)$.  Then
all $(p,q)$ on a given line have the same $u$, $v$, and a representative
is written $(p_u,q_v)$ in the left two columns.   Almost all lines are determined
by their datum $\{[p]_u,[q]_v\}$, with $[ \cdot ]$ indicating reduction modulo
$8$.  The only exceptions are $\{[p]_u,[q]_v\} = \{5_4,1_4\}$ and $\{[p]_u,[q]_v\} = \{1_4,1_4\}$ 
which have two lines each.   The column headed by $\#$ gives the number of
occurrences in our setting $p,q < 250$.  In the five cases where this
number is less than $10$ we continued the computation up through
$p,q < 500$ assuming GRH.  We expect that all possibilities 
are accounted for by the table, and they occur with 
asymptotic frequencies given in the column headed by $\nu$.  
Assuming these frequencies are correct, the average size of $\cK(\mbox{NilOct},-^*p^*q^*)$ is 
exactly $15.875$, substantially larger than the observed $12.01$ in the
$p,q < 250$ setting.  

The connection with  \cite{boston-ellenberg,boston-perry} is as follows.  Let $L(p,q)_k \subset \C$ be the splitting field of 
all degree $2^k$ fields with Galois group a $2$-group and discriminant
$-^* p^* q^*$.    The Galois group $\Gal(L(p,q)_k/\Q)$ is a $2$-group 
and so all ramification at the odd primes $p$ and $q$ is tame.  
Let $L(p,q)$ be the union of these $L(p,q)_k$.   The group $\Gal(L(p,q)/\Q)$ is a pro-$2$-group
generated by the tame ramification elements 
$\tau_p$ and $\tau_q$.    The central question pursued in 
\cite{boston-ellenberg,boston-perry} is the distribution of the $\Gal(L(p,q)/\Q)$ as abstract groups.

Table~\ref{pqoctictab} corresponds to working at the level of the  quotient
$\Gal(L(p,q)_3/\Q)$.  The fact that this group has just the two generators
$\tau_p$ and $\tau_q$ explains why only the sixteen $2$-groups 
having $1$ or $2$ generators appear.  One has $|\Gal(L(p,q)_3/\Q)| = 2^{s_3}$
where $s_3$ is  as in Table~\ref{pqoctictab}.   The lines with an
entry under $T$ are pursued theoretically in 
\cite{boston-ellenberg}.  The cases marked by $\circ$-$\bullet$,
$i$-$iii$, and $I$-$IV$ are respectively treated in 
\S5.2, \S5.3, and $\S5.4$ there.    The entire group 
$\Gal(L(p,q),\Q)$ has order $2^s$,
with $s = \infty$ being expected sometimes in  
Case $IV$.    

Some of the behavior for $k>3$ is previewed by 2-parts of 
class groups of octic fields.  
For example, in Case $ii$ all $87$ instances behave the same:
the unique fields in $\cK(8T2,-^4p^3q^7)$, $\cK(8T4,-^4p^4q^6)$, $\cK(8T17, -^4 p^6 q^5)$, and the
two fields in $\cK(8T17,-^4 p^6 q^7)$ all have $2$ exactly dividing
the class number; the remaining six fields all have odd class number.   
In contrast, in Case $iii$ the $86$ instances 
break into two types of behaviors, 
represented by $(p,q) = (19,5)$ and $(p,q) = (11,37)$.  
These patterns on the database reflect the fact \cite[\S5.3]{boston-ellenberg} that 
in Case $ii$ there is just one possibility for $(\Gal(L(p,q)/\Q);\tau_p,\tau_q)$ 
while in Case $iii$ there are two.

\section{Nilpotent octic fields with discriminant $-^s 2^a q^b$}
\label{octicwild}
The database has all 
octic fields with Galois group 
a $2$-group and discriminant of the form $-^s 2^a q^b$ with 
$q < 2500$.    The sets $\cK(\mbox{NilOct},-^*2^*q^*)$ average
$1711$ fields,
the great increase from the previous section 
being due to the fact that now there are many possibilities for 
wild ramification at $2$.   
As in the previous section, there is great regularity explained by
identifications of relevant absolute Galois groups
\cite{koch}.   Again, even more so this time, there is further regularity not 
explained by theoretical results.   

Continuing with the notation of the previous section, consider the
Galois extensions $L(2,q) = \cup_{k=1}^\infty L(2,q)_k$ and their
associated Galois groups $\Gal(L(2,q)/\Q) = {\varprojlim} \;
\Gal(L(2,q)_k/\Q)$.  As before, octic fields with Galois group a
$2$-group let one study $\Gal(L(2,q)_3/\Q)$.  Table~\ref{2qoctictab}
presents summarizing data for $q < 2500$  in a format parallel to
Table~\ref{pqoctictab} but more condensed.  Here the main entries
count Galois extensions of $\Q$.  Thus an entry $m$ in the $19^2 \;
20 \; 21$ column corresponds to $m$ Galois extensions of $\Q$ having
degree $32$.  Each of these Galois extensions corresponds to four
fields on our database, of types $8T19$, $8T19$, $8T20$, and $8T21$.
 
In the range studied, there are thirteen different behaviors in terms
of the cardinalities $|\cK(8Tj,-^* 2^* q^*)|$.  As indicated by
Table~\ref{pqoctictab}, these cardinalities depend mainly on the
reduction of $q$ modulo $16$.  However classes $1$, $9$, and $15$ are
broken into subclasses.  The biggest subclasses have size $|1A| = 23$, 
$|9A| = 24$, and $|15A|=28$.  The remaining subclasses are 
\begin{align*} 
    1B & =  \{113, 337, 353, 593, 881,1249,1777, 2113, 2129, 2273\}, 
\\
  1C & =  \{257, 1601\}, 
 \\
  1D & = \{577, 1201, 1217, 1553, 1889\}, 
\\
  1E & = \{1153\}, 
\\
  9B & = \{73, 281, 617, 1033, 1049, 1289, 1753, 1801, 1913, 2281,
  2393\},
\\
9C & = \{137, 409, 809, 1129, 1321, 1657, 1993, 2137\}, 
\\
 15B & = \{31, 191, 383, 607, 719, 863, 911, 991,
1103, 1231, 1327, 1471,    \\
&\qquad 1487, 1567, 1583, 2063,  2111, 2287, 2351, 2383\}. 
\end{align*}
A prime $q \equiv 1 \; (16)$ is in $1A$ if and only
if $2 \not \in \F_q^{\times 4}$.  Otherwise we do not have a
concise description of these decompositions.

 \begin{table}[htb]
{\small  
{\renewcommand{\arraycolsep}{1.7pt}
\[
\begin{array}{|c|rrrrr|rrrrrr|rrrrr|rrrr|r|r|}
\hline
  & \multicolumn{5}{c|}{|G|=8} &\multicolumn{6}{c|}{|G|=16} & \multicolumn{4}{c}{|G|=32} &21& \multicolumn{4}{c|}{|G|=64}&&\\
 &  &&&&&&&&&&&&&&&20&&28^2&31^2&&&\\
q & 1 & 2 & 3 & 4 & 5 & 6^2 & 7 & 8 & 9^4 & 10^2 & 11^3 & 
15^2 &16^2 & 17^2 & 18^8 & 19^2 & 26^4 & 27^2 & 29^6 & 30^4 &  35^8 & \mbox{Tot}  \\
\hline
\Q_2 & 24 & 18 & 1 & 18 & 6 & 16 & 36 & 36 & 9 & 12 & 16 & 38 & 12 & 48 & 4 & 24 & 24 & 48 & 16 & 24 & 48 & 1449 \\
\hline \relax
1A &  24 & 18 & 1 & 30 & 2 & 42 & 36 & 44 & 15 & 36 & 12 & 64 & 36 & 96 & 16 & 48 & 80 & 104 & 32 & 52 & 72 & 2895 \\
1B &  \q & \q & \q & \q & \q & 54 & 36 & 60 & \q & \q & \q & 84 & 60 & 144 & \q &
                     96 & 144 & 256 & 80 & 128 & 312 & 6783 \\
1C &  \q& \q & \q & \q & \q & \q & \q & \q & \q & \q& \q& \q& \q & \q & \q &
                     \q & \q & 272 & \q & 136 & 336 & 7071 \\
1D &  \q& \q & \q & \q & \q & \q & \q & \q & \q & \q& \q& \q& \q & \q & \q &
                     \q & \q & 336 & \q & 168 & 384 & 7839 \\
1E &  \q& \q & \q & \q & \q & \q & \q & \q & \q & \q& \q& \q& \q & \q & \q &
                     \q & \q & \q & \q & \q & 240 & 6687 \\
 \hline
9A & 24 & 18 & 1 & 30 & 2 & 44 & 36 & 48 & 15 & 36 & 12 & 68 & 36 & 112 & 16 &
                     48 & 96 & 104 & 48 & 52 & 156 & 3807 \\
9B & \q & \q & \q & \q & \q & \q & \q & \q & \q & \q& \q & \q & \q & \q & \q &
                    \q & \q& \q & \q& \q& 132 & 3615 \\ 
 9C& \q & \q & \q & \q & \q & \q & \q & \q & \q & \q& \q & \q & \q & \q & \q &
                    \q & \q& 72 & \q & 36& 156 & 3615 \\ 
\hline
3,11& 4 & 6 & 1 & 14 & 2 & 10 & 6 & 22 & 7 & 4 & 6 & 21 & 4 & 8 & 3 & 4 & 16 & 8  & 8 & 4 & 21 & 579 \\
\hline
5,13&     8 & 18 & 1 & 12 & 0 & 10 & 20 & 10 & 6 & 12 & 6 & 21 & 12 & 36 & 1 & 12 &   6 & 24 & 4 & 2 & 9 & 621 \\
\hline
7 &  4 & 6 & 1 & 20 & 0 & 30 & 6 & 16 & 10 & 12 & 4 & 34 & 12 & 24 & 8 & 12 & 44 & 24 & 20 & 12 & 60 & 1401  \\
\hline
15A &  4 & 6 & 1 & 20 & 0 & 32 & 6 & 16 & 10 & 12 & 4 & 36 & 16 & 24 & 8 & 20 & 52 & 64 & 24 & 32 & 96 & 2041  \\
15B & \q & \q& \q & \q & \q & \q & \q & \q & \q & \q & \q & \q & \q&\q & \q & \q & \q & \q & \q & \q & 84 & 1945 \\
\hline
\end{array}
\]
}
}
\caption{\label{2qoctictab} The $q$-$j$ entry gives the number of
  Galois extensions of $\Q$ with Galois group $8Tj$ and 
discriminant of the form $-^s 2^a q^b$.  The number of Galois
extensions of $\Q_2$ with 
Galois group $8Tj$ is also given.}
\end{table}

Let $D_{\infty} = \{1,c\}$ where $c$ is complex conjugation.  Let $D_q
\subseteq \Gal(L(2,q)/\Q)$ be a $q$-decomposition group. 
Then, working always in the category of pro-$2$-groups, one has the
presentation  
$D_q = \langle \tau,\sigma | \sigma^{-1} \tau \sigma = \tau^q\rangle$;
here $\tau$ is a ramification 
element and $\sigma$ is a Frobenius element.    Representing a more
general theory, for  
$q \equiv 3,5 \; (8)$ one has two remarkable facts \cite[Example~11.18]{koch}.  First, the $2$-decomposition group $D_2$  
is all of $\Gal(L(2,q)/\Q)$.  Second, the global Galois group is a
free product:  
\begin{equation}
\Gal(L(2,q)/\Q) = D_\infty \ast D_q.
\end{equation}
As a consequence, always for $q \equiv 3,5 \; (8)$, the quotients $\Gal(L(2,q)_k/\Q)$ are computable
as abstract finite groups and moreover depend only on $q$ modulo $8$.   In particular, 
the counts on the lines 3,11 and 5,13 of Table~\ref{2qoctictab} can be obtained purely group-theoretically.
The other lines of Table~\ref{2qoctictab} are not covered by the theory in \cite{koch}.  

A important aspect of the situation is not understood theoretically, namely
the wild ramification at $2$.   The database exhibits extraordinary regularity at the level
$k=3$ as follows.    By $2$-adically completing octic number fields 
$K \in \cK(\mbox{NilOct},-^*2^*q^*)$, one gets $579$ octic $2$-adic fields if
$q \equiv 3 \; (8)$ and $621$ octic $2$-adic fields if $q \equiv 5 \; (2)$.   
The regularity is that the subset of all $1499$ nilpotent octic
$2$-adic fields which arise  
depends on $q$ only modulo $8$, at least in our range $q<2500$. 
 One can see some of this statement
directly from the database: the cardinalities $|\cK(8Tj,-^*2^aq^*)|$ for given
$(j,a)$ depend only on $q$ modulo $8$.  

In the cases $q \equiv 3,5  \; (8)$, the 
group $\Gal(L(2,q)/\Q)=D_2$ has a filtration by higher ramification groups.
From the group-theoretical description of $\Gal(L(2,q)/\Q)$, one can
calculate that 
the quotient group $\Gal(L(2,q)_3/\Q)$ has size $2^{18}$.  The $18$ slopes measuring wildness
of $2$-adic ramification work out to 
\[
\begin{array}{llllll}
3:  \! & 0, & 2, \; 2, \; 2\frac{1}{2} & 3, \; 3, \;  3\frac{1}{2}, \;  3 \frac{1}{2}, \;  3\frac{5}{8}, \;  3 \frac{3}{4}, \; 
& 4, \; 4, \; 4\frac{1}{4}, \; 4 \frac{1}{4}, \; 4 \frac{3}{8}, \; 4\frac{1}{2},   \;  4\frac{3}{4} & 5 \\
5:  \! & 0, \; 0, \; & 2, \; 2, \; 2, \; 2\frac{1}{2}& 3, \; 3, \;  3, \;  3 \frac{1}{2}, \;  3\frac{1}{2}, \;  3 \frac{3}{4}, \;
&  4, \; 4\frac{1}{4}, \; 4\frac{1}{2},   \;  4\frac{3}{4}, \;  4 \frac{3}{4}, \; & 5. 
\end{array}
\]
Most of these slopes can be read off from the octic field part of the
 database directly, via the automatic $2$-adic analysis of fields given there.
  For example, the first four slopes 
for $q=3$ all arise already from $\Q[x]/(x^8+6x^4-3)$, the unique member of $\cK(8T8,-^3 2^{16} 3^7)$.
A few of the listed slopes can only be seen directly by working with degree sixteen resolvents.   
A natural question, not addressed in the literature, is to similarly describe the 
slopes appearing in all of $\Gal(L(2,q)/\Q)$.

\section{Minimal nonsolvable fields with  $\grd \leq \Omega$} 
\label{minnonsolvable}
   Our focus for the remainder of the paper is on Galois number fields, for
which root discriminants and Galois root discriminants naturally 
coincide.  As reviewed in the introduction, in \cite{jrlowgrd} we raised 
the problem of completely understanding the set $\cK[\Omega]$ of
all Galois number fields $K \subset \C$ with $\GRD$ at most the
Serre-Odlyzko constant $\Omega = 8 \pi e^\gamma \approx 44.76$.   
As in \cite{jrlowgrd}, we focus attention here on the interesting
subproblem of identifying the subset $\cK^{\rm ns}[\Omega]$
of $K$ which are nonsolvable.  
Our last two sections explain how the database explicitly exhibits a substantial
part of $\cK^{\rm ns}[\Omega]$.    

 \begin{figure}[htb]
\includegraphics[width=5in]{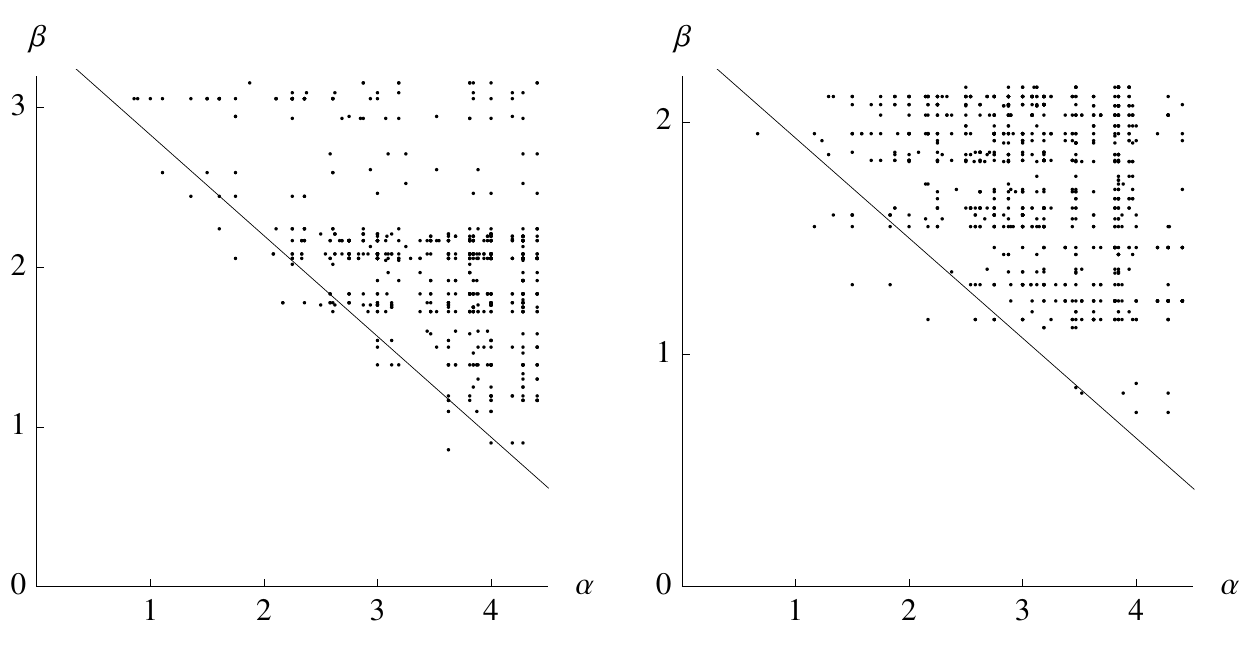}
\caption{\label{nonsolvefig} Galois root discriminants $2^\alpha 3^\beta$ (left) and $2^\alpha 5^\beta$ (right) 
arising from minimal nonsolvable fields of degree $\leq 11$ in the database.  The lines have equation
$2^\alpha q^\beta = \Omega$.} 
\end{figure}

   We say a nonsolvable number field is {\em minimal} if it
does not contain a strictly smaller nonsolvable number 
field.   So fields with Galois group say  $S_n$ are minimal, 
while fields with Galois group say
 $C_p \times S_n$ or $C_p^k:S_n$ are not.  
Figure~\ref{nonsolvefig} draws a dot for each
minimal nonsolvable field $K_1 \in \cK^{\rm ns}_{\rm min}[\Omega]$
coming the degree $\leq 11$ part of the database 
with $\GRD$ of the form $2^\alpha 3^\beta$
or $2^\alpha 5^\beta$.  There are 654 fields in the
first case and 885 in the second.  Of these fields,
24 and 17 have $\GRD \leq \Omega$.  Figure~\ref{nonsolvefig} 
illustrates the extreme extent to which the low $\GRD$ problem is focused 
on the least ramified of all Galois number fields.

    Figure~\ref{nonsolvefig} also provides some context for the next section
    as follows.   Consider the compositum $K = K_1 K_2$ of
 distinct minimal fields $K_1$ and $K_2$ contributing to the same
half of Figure~\ref{nonsolvefig}.  Let $2^{\alpha_i} q^{\beta_i}$ be the 
root discriminant of $K_i$.  The root discriminant $2^\alpha q^\beta$
 of $K$ satisfies $\alpha \geq \max(\alpha_1,\alpha_2)$ and
 $\beta \geq \max(\beta_1,\beta_2)$.    The figure makes it clear that one
 must have almost exact agreement
  $\alpha_1 \approx \alpha_2$ and $\beta_1 \approx \beta_2$ for
 $K$ to even have a chance of lying in $\cK^{\rm ns}[\Omega]$.   
 As some examples where one has exact agreement, consider the respective splitting fields 
 $K_1$, $K_2$, and $K_3$ of
  \begin{eqnarray*}
  f_1(x) & = & x^5-10 x^3-20 x^2+110 x+116, \\
  f_2(x) & = & x^5+10 x^3-10 x^2+35 x-18, \\
  f_3(x) & = &  x^5+10 x^3-40 x^2+60 x-32.
  \end{eqnarray*}
  All three fields have Galois group $S_5$ and root discriminant $2^{3/2} 5^{8/5} \approx 37.14$.  
  The first two completely agree at $2$, but differ at $5$, so that $K_1K_2$ has root discriminant 
  $2^{3/2} 5^{48/25} \approx 62.17$.   The other two composita also have root discriminant
  well over $\Omega$, with $\GRD(K_2 K_3) = 2^{9/4} 5^{8/5} \approx 62.47$ and 
  $\GRD(K_1 K_3) = 2^{9/4} 5^{48/25} \approx 104.55$.  These computations,
  done automatically by entering $f_i(x) f_j(x)$ into the $\GRD$ calculator of
  \cite{jr-local-database}, are clear illustrations of the general
  difficulty of  
  using known fields in $\cK[\Omega]$ to obtain others.  

In
 \cite{jrlowgrd}, we listed fields proving 
 $|\cK_{\rm min}^{\rm ns}[\Omega]| \geq 373$.
 Presently, the fields on the database show
 $|\cK_{\rm min}^{\rm ns}[\Omega]| \geq 386$.  
 In \cite{jrlowgrd}, we highlighted the fact
 that the only simple groups involved were
the five smallest, $A_5$, $SL_3(2)$, $A_6$, $SL_2(8)$
and $PSL_2(11)$,  and the eighth, $A_7$.  
The new fields add $SL_2(16)$
, $G_2(2)'$, and $A_8$ to the 
list of simple groups involved.  These groups are $10^{\rm th}$, $12^{\rm th}$, and 
tied for $19^{\rm th}$ on the list of all non-abelian simple groups in 
increasing order of size.   

\begin{table}[htb]
\[
\begin{array}{|r|r|cr|cr|}
\hline
\# & |H| & G=H & \# & G=H.Q & \# \\
\hline
1 & 60 & A_5 & {\bf 78} & S_5 & {\bf 192} \\
2 & 168 & SL_3(2) & .\mathit{18} &  PGL_2(7) &...\mathit{23}  \\ 
3 & 360 & A_6 &   {\bf 5} & S_6, \; PGL_2(9), \; M_{10}, \; P\Gamma L_2(9) \; & \!\! {\bf 13}, \; ....\mathit{6}, \; \mathit{0}, \; .\mathit{15} \\
4 & 504 & SL_2(8) & \mathit{15} & \Sigma L_2(8) & \mathit{15} \\
5 & 660 & PSL_2(11) & \mathit{1} & PGL_2(11) & \mathit{0}\\
8 & 2520 & A_7 & \mathit{1} &  S_7 & \mathit{1} \\
1^2 & 3600 & A_5^2 & & A_5^2.2, \; A_5^2.V, \;  A_5^2.C_4,  \; A_5^2.D_4 & \mathit{1}, \; .\mathit{1},\;  \mathit{0}, \; \mathit{0}  \\
10 & 4080 & SL_2(16) & .\mathit{1} &  SL_2(16).2, SL_2(16).4 & \mathit{0},\mathit{0} \\
12 & 6048 & G_2(2)' & \mathit{0} & G_2(2) & .\mathit{1} \\
19 & 20160 & A_8 & \mathit{0} & S_8 & .\mathit{1} \\
\hline
\end{array}
\]
\caption{\label{grdsummary} Lower bounds on $|\cK[G,\Omega]|$ for minimal nonsolvable groups $G$.
Entries highlighted in bold are completeness results from \cite{jrlowgrd}.  Fields found since 
\cite{jrlowgrd} are indicated by .'s.  }
\end{table}

     Table~\ref{grdsummary} summarizes all fields on the database in $\cK_{\rm min}^{\rm ns}[\Omega]$.  It is organized by the
      socle $H \subseteq G$, which is a simple group except in the single case $H = A_5 \times A_5$.   The .'s indicate that, for example, 
      of the $23$ known fields in $\cK[PGL_2(7),\Omega]$, twenty
are listed in \cite{jrlowgrd} and three are new.   The polynomial for
the $SL_2(16)$ field was found by Bosman \cite{bosman}, starting from 
a classical modular form of weight $2$.   We found polynomials for the
new $SL_3(2)$ field and the three new $PGL_2(7)$ fields 
starting from Schaeffer's list \cite[App~A]{schaeffer} of ethereal modular forms of weight $1$.  Polynomials for the
other new fields were found by specializing families.  All fields summarized by Table~\ref{grdsummary} come
from 
the part of the database in degree $\leq 11$, except for Bosman's degree seventeen polynomial and 
the degree twenty-eight polynomial for $G_2(2)$.   It would be of interest to pursue
calculations with modular forms more systematically.  They have the potential not only to yield new fields
in $\cK_{\rm min}^{\rm ns}[\Omega]$, but also to prove completeness for certain $G$.

\section{General nonsolvable fields with $\GRD \leq \Omega$}
\label{gennonsolvable}

We continue in the framework of the previous section, so that the focus remains
on Galois number fields contained in $\C$.  
    For $K_1 \in \cK_{\rm min}^{\rm ns}[\Omega]$ such a Galois number field, let $\cK[K_1;\Omega]$ be the
subset of $\cK^{\rm ns}[\Omega]$ consisting of fields containing $K_1$.  
Clearly 
\begin{equation}
\label{subodunion}
\cK^{\rm ns}[\Omega] = \bigcup_{K_1} \cK[K_1;\Omega].
\end{equation}
So a natural approach to studying all of $\cK^{\rm ns}[\Omega]$ is to study 
each $\cK[K_1;\Omega]$ separately.

The refined local information contained in the database 
can be used to find fields in 
$\cK[K_1;\Omega]$.   The set of fields
so obtained is  always very small, often just 
$\{K_1\}$.   Usually it seems likely that the
set of fields obtained is all of $\cK[K_1;\Omega]$,
and sometimes this expectation is provable 
under GRH.   We sketch such a proof for a particular $K_1$ in the first example below.
In the remaining examples, we start from other $K_1$ 
and now construct 
proper extensions $K \in \cK[K_1;\Omega]$,
 illustrating several phenomena. 
Our examples are organized in terms of
increasing degree $[K:\Q]$.   The fields here
are all extremely lightly ramified for their Galois
group, and therefore worthy of individual attention.  

Our local analysis of a Galois number field $K$ centers on the notion
of $p$-adic slope content described in \cite[\S3.4]{jr-local-database}
and automated on the associated database.  Thus a $p$-adic slope
content of $[s_1,\dots,s_m]_t^u$ indicates a wild inertia group $P$ of
order $p^m$, a tame inertia group $I/P$ of order $t$, and an
unramified quotient $D/I$ of order $u$.  Wild slopes $s_i \in \Q \cap
(1,\infty)$ are listed in weakly increasing order and from 
\cite[Eq.~7]{jr-local-database} the contribution $p^\alpha$ to the root
discriminant of $K$ is determined by
\[
\alpha = \left( \sum_{i=1}^m \frac{p-1}{p^i} s_{n+1-i} \right) + \frac{1}{p^m} \frac{t-1}{t}.
\]
The quantities $t$ and $u$ are omitted from presentations of slope content when they are $1$.
  
\subsubsection*{Degree $120$ and nothing more from $S_5$} The polynomial
\[
f_1(x) = x^5 + x^3 + x - 1
\]
has splitting field $K_1$ with root discriminant $\Delta_1 = 11^{2/3}
37^{1/2} \approx 30.09$.  
Since $\Delta_1 2^{2/3} \approx 47.76$, $\Delta_1 3^{1/2} \approx
52.11$, $\Delta_1 11^{1/6} \approx 44.87$, and $\Delta_1 37^{1/4}
\approx 74.20$ are all more than $\Omega$, any $K \in \cK[K_1;\Omega]$
has to have root discriminant $\Delta = \Delta_1$.  The GRH bounds say
that a field with root discriminant $30.09$ can have degree at most
$2400$ \cite{martinet}.

The main part of the argument is to use the database to 
show that most other {\em a priori} possible
$G$ in fact do not arise as $\Gal(K/\Q)$ for 
$K \in \cK[K_1;\Omega]$.     For example, 
if there were an $S_3$ field $K_2$ with absolute
discriminant $11^2 37$, then $K_1K_2$ would be in 
$\cK[K_1; \Omega]$; there is in fact an $S_3$ field
with absolute discriminant $11 \cdot 37^2$, but not one
with absolute discriminant $11^2 37$.      
As an example of a group that needs a supplementary argument
to be eliminated,  
the central extension $G = 2.S_5$ does not appear because the degree 12 subfield of $K_1$ fixed 
by $D_5 \subset S_5$ has root discriminant $\Delta_1$ and 
 class number $1$.  

\subsubsection*{Degree $1920$ from $A_5$}  The 
smallest root discriminant of any nonsolvable Galois field
is $2^{6/7} 17^{2/3} \approx 18.70$ coming
from a field $K_1$ with Galois group $A_5$. 
This case is complicated because one can add ramification 
in several incompatible directions, so that there are different
maximal fields in $\cK[K_1;\Omega]$.  One overfield is 
the splitting field $\tilde{K}_1$ of $f_-(x)$ where
\[
f_{\pm}(x) = x^{10} + 2 x^6 \pm 4 x^4 - 3 x^2  \pm 4.
\]
In this direction, ramification has been added at $2$ making
the slope content there $[2,2,2,2,4]^6$ and the
root discriminant $2^{39/16}  17^{2/3} \approx 35.81$.  
The only solvable field $K_2$ on the database which
is not contained in $\tilde{K}_1$ but has $\mbox{rd}(\tilde{K}_1 K_2) < \Omega$ 
is $\Q(i)$.   The field $\tilde{K}_1 K_2$ is the splitting field
of $f_+(x)$ with Galois group $10T36$.  There is yet
another wild slope of $2$, making the root discriminant $2^{79/32} 17^{2/3} \approx 36.60$.  

\subsubsection*{Degree $25080$ from $PSL_2(11)$.}  The only known field $K_1$ with 
Galois group $PSL_2(11)$ and root discriminant less than $\Omega$ first appeared in \cite{kluners-malle} and
 is the splitting field
of 
\[
f_1(x) = x^{11}-2 x^{10}+3 x^9+2 x^8-5 x^7+16 x^6-10 x^5 +10 x^4+2 x^3-3 x^2+4 x-1.
\]
The root discriminant is $\Delta_1 = 1831^{1/2} \approx 42.79$, forcing all members of 
$\cK[K_1; \Omega]$ to have root discriminant $1831^{1/2}$ as well.  

The prime $1831$ is congruent to $3$ modulo $4$, so that the associated
quadratic field $\Q(\sqrt{-1831})$ is imaginary and its class number 
can be expected to be considerably larger than $1$.  This class number is in fact $19$, and 
the splitting field of a degree $19$ polynomial in the database is the
corresponding Hilbert class field $K_2$.  The field $K_1 K_2 \in \cK[K_1;\Omega]$ has degree $660 \cdot 38 = 25080$.

\subsubsection*{Degree $48384$ from $SL_2(8).3$.}  The splitting field $K_1$ of
\[
f_1(x) = x^9-3 x^8+4 x^7+16 x^2+8 x+8 
\]
has Galois group $\Gal(K_1/\Q) = 9T32 = SL_2(8).3$ and root discriminant $2^{73/28} 7^{8/9} \approx 34.36$.
This root discriminant is the smallest known from a field with Galois group $SL_2(8).3$.  In fact,
it is small enough that it is possible to add ramification at both $2$ and $7$ and still keep 
the root discriminant less than $\Omega$.  Namely let 
\begin{eqnarray*}
f_{2}(x) & = & x^4-2 x^3+2 x^2+2, \\
f_{3}(x) & = & x^4-x^3+3 x^2-4 x+2.
\end{eqnarray*}
The splitting fields $K_{2}$ and $K_{3}$ have Galois groups
$A_4$ and $D_4$ respectively.     Composing with $K_{2}$
increases degree by four and adds wild slopes $2$ and $2$ to the original
$2$-adic slope content $[20/7,20/7,20/7]^3_7$.   Composing with $K_{3}$ 
then increases degrees by $8$, adding another wild slope of $2$ to the $2$-adic slope content
and increasing the $7$-adic tame degree from $9$ to $36$.  The root discriminant of $K_1 K_{2} K_{3}$ 
is then $2^{153/56} 7^{35/36} \approx 44.06$.  

\subsubsection*{Degree $80640$ from $S_8$}   The largest group in Table~\ref{grdsummary} is $S_8$, and
the only known field in $\cK[S_8,\Omega]$ is the splitting field $K_1$ of
\[
f_1(x) = x^8-4 x^7+4 x^6+8 x^3-32 x^2+32 x-20.
\]
Here Galois slope contents are $[15/4,7/2,7/2,3,2,2]^3$ and $[ \; ]_7$ at $2$ and $5$ respectively,
giving root discriminant $2^{111/32} 5^{6/7} \approx 43.99$.  The only field on the database which can be used
to give a larger field in $\cK[K_1; \Omega]$ is $K_2 = \Q(i)$.  This field gives an extra wild slope of
$2$, raising the degree of $K_1 K_2$ to $80640$ and the root discriminant to $2^{223/64} 5^{6/7} \approx 44.47$.  

\subsubsection*{Degree $86400$ from $A_5^2.V$} 
Another new field $K_1$ on Table~\ref{grdsummary}, found by Driver, 
is the splitting field of
\[
f_1(x) = x^{10}-2 x^9+5 x^8-10 x^6+28 x^5-26 x^4 -5 x^2+50 x-25.
\]
Like in the previous example, this field $K_1$ is wildly ramified at $2$ and tamely ramified at $5$.
Slope contents are $[23/6,23/6,3,8/3,8/3]_3$ and $[ \; ]_6$ for a root discriminant
of $2^{169/48} 5^{5/6} \approx 43.89$.  The splitting field $K_2$ of $x^3-x^2+2 x+2$ 
has Galois group $S_3$, with $2$-adic slope content $[3]$ and $5$-adic slope content 
$[\;]_3$.   In the compositum $K_1 K_2$, the extra slope is in fact $2$ giving a root discriminant of
$2^{85/24} 5^{5/6} \approx 44.53$.  

\subsubsection*{Degree $172800$ from $S_5$ and $S_6$}  Consider the $\binom{386}{2} = 74305$
composita $K_1 K_2$, as $K_1$ and $K_2$ vary over distinct known fields in 
$\cK^{\rm ns}_{\rm min}[\Omega]$.   From our discussion of Figure~\ref{nonsolvefig}, one 
would expect that very few of these  composita would have root discriminant less than 
$\Omega$.  In fact, calculation shows that exactly one of these composita has $\mbox{rd}(K_1 K_2) \leq \Omega$, namely the
joint splitting field of
\begin{eqnarray*}
f_1(x) & = & x^5 - x^4 - x^3 + 3 x^2 - x - 19, \\
f_2(x) & = & x^6 - 2 x^5 + 4 x^4 - 8 x^3 + 2 x^2 + 24 x - 20.
\end{eqnarray*}
Here $\Gal(K_1/\Q) = S_5$ and $\Gal(K_2/\Q) = S_6$.  Both fields have tame ramification of order $2$ at $3$ and
order $5$ at $7$.  Both are otherwise ramified only at $2$, with $K_1$ having slope content $[2,3]^2$ and
$K_2$ having slope content $[2,2,3]^3$.   In the compositum $K_1 K_2$, there is partial cancellation between the two wild slopes 
of $3$, and the slope content is $[2,2,2,2,3]^6$.   The root discriminant of $K_1 K_2$ then works out to 
$2^{39/16} 3^{1/2} 7^{4/5} \approx 44.50$.  The existence of this remarkable compositum
contradicts Corollary~12.1 of \cite{jrlowgrd} and is the only error we have found in
\cite{jrlowgrd}.  

The field discriminants of $f_1$ and $f_2$ are respectively $-^2 2^6 3^1 7^4$ and $-^2 2^9 3^1 7^4$.   
The splitting fields $K_1$ and $K_2$ thus contain distinct quadratic fields,  $\Q(\sqrt{3})$ and $\Q(\sqrt{6})$ respectively.  
The compositum therefore has Galois group all of $S_5 \times S_6$, and so the degree $[K_1K_2 : \Q] = 120 \cdot 720 = 86400$
ties that of the previous example.  But, moreover,  $K_3 = \Q(\sqrt{-3})$ is disjoint from $\Q(\sqrt{3},\sqrt{6})$ and
does not introduce more ramification.  So $K = K_1 K_2 K_3$ has the same root discriminant $2^{39/16} 3^{1/2} 7^{4/5} \approx 44.50$, but the larger degree $2 \cdot 86400 = 172800$.

\medskip

The GRH upper bound on degree for a given root discriminant $\delta \in [1,\Omega)$
increases to infinity as $\delta$ increases to $\Omega$ (as illustrated by Figure 4.1 of \cite{jrlowgrd}).  
 However, we have only
exhibited fields $K$ here of degree $\leq 172800$.  Dropping the restriction that
$K$ is Galois and nonsolvable may let one  obtain somewhat larger degrees,
 but  
there remains a substantial and intriguing gap between degrees of known fields and analytic upper
bounds on degree.

\bibliographystyle{amsalpha}
\bibliography{jr}

\end{document}